\documentclass[11pt,a4paper]{amsart}
 \pdfoutput=1
\usepackage{graphicx}
\usepackage{amsmath}
\usepackage{algorithmicx}
\usepackage{algpseudocode}
\usepackage{algorithm}

\usepackage{pgfplots,tikz}
\usetikzlibrary{shapes,arrows,calc}
\usepackage[centering]{geometry}
\usepackage{etoolbox}

\usepackage{natbib}

\usepackage{tikz}
\usetikzlibrary{matrix,shapes,calc,through,positioning,chains,automata,arrows,intersections,decorations.pathmorphing,backgrounds,fit,fadings,decorations.pathreplacing,patterns,shadows,trees,decorations.markings}
\tikzstyle{arrow} = [draw, -latex']
\tikzstyle{every edge}=  [draw=red,fill=blue!20,opacity=.5]  
\tikzset{
	block/.style=
	{
		draw, rectangle, minimum height=0.9cm, minimum width=1cm
	},
	bblock/.style=
	{
		draw, rectangle,minimum height=1cm, minimum width=1.4cm
	},
	bbblock/.style=
	{
		draw, rectangle,minimum height=1.4cm, minimum width=1.6cm 
	},
	bbbblock/.style=
	{
		draw, rectangle,minimum height=2cm, minimum width=2.2cm 
	},
	block_nl/.style=
	{
		draw,double, rectangle, minimum height=0.9cm, minimum width=1cm
	},
	bblock_nl/.style=
	{
		draw,double, rectangle,minimum height=1cm, minimum width=1.4cm
	},
	bbblock_nl/.style=	
	{
		draw,double, rectangle,minimum height=1.1cm, minimum width=1.6cm 
	},
	branch/.style=
	{
		circle, inner sep=0pt, minimum size=1mm, fill=black, draw=black
	},
	connector/.style=
	{
		->, thick
	},
	connector2/.style=
	{
		->, dashed
	},
	dummy/.style=
	{
		inner sep=0pt, minimum size=0pt
	},
	dummyfilled/.style=
	{
		inner sep=0pt, minimum size=0.4pt, fill=black, draw=black
	},
	inout/.style=
	{
	},
	sum/.style=
	{
		circle, inner sep=0pt, minimum size=2mm, draw=black, thick
	},
	sample1/.style=
	{
		circle, inner sep=0pt, minimum size=1mm, draw=black, thick
	},
	sample2/.style=
	{
		circle, inner sep=0pt, minimum size=1.2cm, draw=black, thick
	}, 
	vecArrow/.style =
	{
		thick, decoration={markings,mark=at position
			1 with {\arrow[semithick]{open triangle 60}}},
		double distance=1.4pt, shorten >= 5.5pt,
		preaction = {decorate},
		postaction = {draw,line width=1.4pt, white,shorten >= 4.5pt}
	},
	innerWhite/.style =
	{semithick, white,line width=1.4pt, shorten >= 4.5pt
	},
	kasten/.style=
	{
		rectangle,draw=black!90,fill=black!10,thick,inner sep=2pt,minimum size=7mm
	},
	kreis/.style=
	{
		circle,draw=black!90,fill=black!10,thick,inner sep=0pt,minimum size=2mm
	},
	point/.style=
	{
		{coordinate},very thick,draw=black
	},
	kreuz/.style=
	{
		circle,draw=black,fill=black,thick,inner sep=0pt,minimum size=4pt
	},
	kreuz_inv/.style=
	{
		circle,draw=white,fill=white,thick,inner sep=0pt,minimum size=4pt
	},
}

\usepackage{stmaryrd}
\usepackage{color}
\definecolor{darkgreen}{rgb}{0.0, 0.5, 0.0}
\def\LG{\textcolor{blue}}

\newcommand{\mbb}[1]{\mathbb #1}

\newcommand{\mcl}[1]{\mathcal #1}

\newtheorem{definition}{Definition}
\newtheorem{remark}{Remark}

\makeatletter
\newcommand{\fixed@sra}{$\vrule height 2\fontdimen22\textfont2 width 0pt\shortrightarrow$}
\newcommand{\shortarrow}[1]{%
  \mathrel{\text{\rotatebox[origin=c]{\numexpr#1*45}{\fixed@sra}}}
}
\makeatother
\newcommand{\boundellipse}[3]
{(#1) ellipse (#2 and #3)
}

\begin{document}

\title[Turnpike Properties in Optimal Control]{Turnpike Properties in Optimal Control\\
{\small An overview of discrete-time and continuous-time results}}
\author[T.\ Faulwasser \& L. Gr\"u{}ne]{Timm Faulwasser$^{1}$ and Lars Gr\"u{}ne$^{2}$}
	\thanks{}
	\thanks{$^{1}$TU Dortmund University, Institute of Energy Systems, Energy Efficiency and Energy Economics, Germany
		{\tt\small timm.faulwasser@ieee.org}}%
	\thanks{$^{2}$Chair of Applied Mathematics, Mathematical Institute, Universität Bayreuth, 
   Bayreuth, Germany {\tt\small lars.gruene@uni-bayreuth.de}}
	\maketitle

\begin{abstract}
	The turnpike property refers to the phenomenon that in many optimal control problems, the solutions for different initial conditions and varying horizons approach a neighborhood of a specific steady state, then stay in this neighborhood for the major part of the time horizon, until they may finally depart. While early observations of the phenomenon can be traced back to works of Ramsey and von Neumann on problems in economics in 1928 and 1938, the turnpike property received continuous interest in economics since the 1960s and recent interest in systems and control. The present chapter provides an introductory overview of discrete-time and continuous-time results in finite and infinite-dimensions. We comment on dissipativity-based approaches and infinite-horizon results, which enable the exploitation of turnpike properties for the numerical solution of problems with long and infinite horizons. After drawing upon numerical examples, the chapter concludes with an outlook on time-varying,  discounted, and open problems.  
\end{abstract}
\smallskip
\noindent \textit{Keywords.} Turnpike properties, Optimal control, Dissipativity, Numerical solution

\section{Introduction and Historical Origins}\label{sec:Intro}
Optimal control  has a rich history, which can be traced back at least to the 1600s \citep{Sussmann97}. Moreover, around the mid of the 20th century two complementary breakthroughs---the Pontryagin Minimum Principle (PMP)~\citep{Boltyanskii60a} and Dynamic Programming (DP)~\citep{Bellman54a}---have been achieved for which the investigations into the Calculus of Variations have paved the road, see \citep{Pesch13a}. The number of papers, books, and PhD theses triggered by the PMP and DP is countless. 
This evident scientific appeal  is likely due to optimal control being a unique combination of powerful mathematical concepts applicable to the entire spectrum of dynamics systems from finite to infinite dimensions---including differential algebraic settings, hybrid formulations, and inclusions. At the same time, optimal control is nowadays readily applicable in a variety of contexts ranging from analysis of phenomena in nature and society towards the design and automation of technical systems. Hence, it comes without surprise that the list of successful applications of optimal control is ever growing since and includes highlights such as: the moon landing in 1969 and numerous space missions; insights into  brain evolution in humans~\citep{Gonzalez18a}; understanding the complex interplay between economic growth, climate change, and green house gas emissions~\citep{Nordhaus-1992}; design and understanding of methods for machine learning~\citep{EHJ17,EGPZ20,Li17a};  as well as numerous mass-consumer products, which entail feedback control loops designed and implemented via optimal control theory. In the view of automatic control applications, a tremendously successful off-spring of the family of methods subsumed as optimal control is its receding horizon variant, which is also denoted as model predictive control~\citep{Rawlings17}. 

\subsection*{From Ramsey and von Neumann to Aircraft Design}

This chapter is concerned with providing an overview on the turnpike phenomenon and its numerical application, which occurs frequently in real-world Optimal Control Problems~(OCPs). It refers to qualitative and quantitative similarity of optimal solutions to OCPs for different initial
conditions and varying horizons. Specifically, the solutions approach a neighbourhood of the best steady
state, stay within this neighbourhood for some time, and
might leave this neighbourhood towards the end of the
optimization horizon. The notion of turnpike properties in OCPs has been coined by \cite[p. 331]{Dorfman58} with the following lines:

\begin{center}
\begin{minipage}[h]{0.95\textwidth}
\noindent\textit{It is exactly like a turnpike paralleled by a network of minor roads. [...], if origin and destination are far enough part, it will always pay to get on to the turnpike and cover distance at the best rate of travel, even if this means adding a little mileage at either end}.
\end{minipage}
\end{center}

Prior to this foundational description, the phenomenon has been reported by  \cite{vonNeumann38} and by \cite{Ramsey28}. Indeed, without using the term \textit{turnpike}, \cite{Samuelson49} implicitly conjectured corresponding theorems. Given that these early works on the subject investigated specific problems arising in economics, it is not surprising that turnpike theory was driven by research on applying optimal control in economics and econometrics~\citep{Mckenzie76,CaHL91}. Interestingly, also in systems and control the phenomenon was observed at least as early as~\citep{Wilde72a}---though Kokotovic and Wilde used a different name, referring to it as a \textit{dichotomy in optimal control}.\footnote{More recently, the name \textit{hyper-sensitive} optimal control problems was also used \citep{Rao99a}.}
Historically, the majority of turnpike results look at specific problems and exploit structure of the optimality conditions. Recent papers have proposed to leverage turnpikes in applications as diverse as shape optimization in aircraft design~\citep{LaTZ20}, membrane-filtration systems~\citep{Kalboussi18}, or control of chemical reactors with uncertain models~\citep{kit:faulwasser19c}. 

Remarkably, in the systems and control community, turnpikes have received little attention. One of the earliest comment on turnpike properties in receding-horizon optimal control appears to be~\citep{Rawlings08a}, who later coined the phrase economic MPC~\citep{Rawlings09b} for receding horizon optimal control with generic objectives.\footnote{\textit{Generic} refers to objectives that do not merely penalize the distance to some reference setpoint.} Eventually, the interest in economic MPC has led to a new line of research towards turnpike properties centred around dissipativity notions, which have been introduced by Jan Willems~\citep{Willems71a,Will72a,Will72b}. 
Early works in this direction include~\citep{Grue13,epfl:faulwasser14e}, while \cite{GruM16} and \cite{FKJB17} investigate equivalence of turnpike and dissipativitiy properties of OCPs in discrete-time, respectively, continuous-time settings. Given the close relation between dissipativity and stability of dynamic systems, it is far from surprise, that the recent elaboration of the dissipativity route to turnpikes establishes a close relation between turnpikes and stability in infinite-horizon optimal control~\citep{GrKW17,tudo:faulwasser20a}. All these works refer to steady-state turnpikes. We will comment on the progress with respect to time-varying and periodic turnpikes in the outlook in Section~\ref{sec:open_problems}.

The remainder of this chapter is structured as follows: In Section \ref{sec:Intro} we conclude with the problem setup and a review of the necessary optimality conditions.  Section~\ref{sec:Taxonomy} attempts a taxonomy of turnpike properties, this way providing guidance in the considerable amount of literature on this subject.  Section~\ref{sec:Results} details the generating mechanisms as well as necessary and sufficient conditions for the occurrence of the phenomenon. Finally, Section~\ref{sec:Numerics} comments on the exploitation of turnpikes in numerics and receding-horizon optimal control. This chapter concludes with an outlook on advanced turnpike settings and open problems of interest.

\subsection*{Notation}
The set  $\mathbb{N}_{[a,b]} \doteq \{a,\ldots,b\} \subset \mbb{N}$ denotes an interval of natural numbers.  
Functions $u:[0,T] \to \mcl{U}$ and sequences $u:\{a, a+1, \dots, b-1, b\} \to \mcl{U}$ are written as $u(\cdot)$. Partial derivatives are indicated by subscripts, i.e. $\ell_x$ is the partial derivative of $\ell$ w.r.t. $x$. Throughout this chapter, optimal solutions to optimization and optimal control problems are indicated by superscript $(\cdot)^\star$.  Variables at steady state are marked by superscript $\bar{(\cdot)}$; optimal steady-state variables are indicated by $(\bar{\cdot})^\star$.
 
\subsection{Discrete-time and Continuous-time Setting }

We consider discrete-time OCPs of the following form 
\begin{subequations}
	\label{eq:ocp_dt}
	\begin{align}
V_T(x_0) \doteq\min_{u(\cdot)} \ \ & \sum_{t=0}^{T-1} \ell(x(t),u(t))  + \varphi(x(T))  \label{eq:ocp_dt_obj} \\
	\text{subject to} \nonumber \\  
	 x(t+1) &= f(x(t),u(t)), ~x(0) = x_0&&\hspace{-0.5em} t\in \mathbb{N}_{[0,T-1]}, \label{eq:ocp_dt_dyn}\\
	 0&\geq g(x(t),u(t)), &&\hspace{-0.5em}  t\in \mathbb{N}_{[0,T-1]}, \label{eq:ocp_dt_pc}\\
	 0&\geq \psi(x(T)). \label{eq:ocp_dt_tc}
	\end{align}
\end{subequations}
The constraint set for $z\doteq (x, u) \in\mcl{X} \times \mcl{U}$, where $\mcl{X}$ and $\mcl{U}$ are normed metric spaces,  is defined as
\begin{equation} \label{eq:Zset}
\mbb{Z} = \{ z= (x, u) \in\mcl{X} \times \mcl{U}\,|\, g_j(z) \leq 0, \, j \in \mathbb{N}_{[1,n_g]}\},
\end{equation}
where the constraints \eqref{eq:ocp_dt_pc} are given by $ g(x,u) = \begin{bmatrix}
g_1(x,u) & \dots &g_{n_g}(x,u)\end{bmatrix}^\top$ with $g_j:\mcl{X} \times \mcl{U} \to \mbb{R}$,  $j \in \mathbb{N}_{[1,n_g]}$.

The projection of $\mbb{Z}$ onto $\mcl{X}$ is denoted as $\mbb{X}$ and the projection of $\mbb{Z}$ onto $\mcl{U}$ is written as $\mbb{U}$.\footnote{In case this projection depends on the state $x$, we write $\mbb{U}(x)$. }  The continuous functions $\ell:\mcl{X}\times\mcl{U}\to\mbb{R}$ and $\varphi: \mcl{X}\to\mbb{R}$ are the Lagrange function and, respectively,  the Mayer term. Moreover, $T \in \mbb{N}$ is the horizon. 

We remark that, throughout this chapter, unless otherwise stated, we do not make specific assumptions on $\ell$ and $\varphi$ such as quadratic  structures or convexity. Instead we will work with a rather generic dissipativity assumption to be introduced later. Equation \eqref{eq:ocp_dt_tc} defines the terminal constraint; the more general case of coupled constraints on initial and terminal states is left out due to space limitations. 

The continuous-time counterpart to OCP \eqref{eq:ocp_dt} reads
\begin{subequations}
	\label{eq:ocp_ct}
	\begin{align}
V_T(x_0) \doteq\min_{u(\cdot) \in \mcl{L}([0,T], \mcl{U})} \ \ & \int_{0}^{T} \ell(x(t),u(t))\mathrm{d}t  + \varphi(x(T))  \label{eq:ocp_ct_obj} \\
	\text{subject to} \nonumber \\  
	 \dot x &= f(x,u), ~x(0) = x_0&&\hspace{-0.5em} t \in [0,T], \label{eq:ocp_ct_dyn}\\
	 0&\geq g(x(t), u(t)), &&\hspace{-0.5em}  t \in [0, T], \label{eq:ocp_ct_pc}\\
	 0&\geq \psi(x(T)) . \label{eq:ocp_ct_tc}
	\end{align}
\end{subequations}
The projected constraint sets are similarly defined as in case of OCP \eqref{eq:ocp_dt}. The horizon is $T\in\mbb{R}^+$. However, the underlying space is considered to be a real-valued vector space of appropriate finite dimension, i.e., \eqref{eq:ocp_ct} refers to a finite-dimensional OCP. In order guarantee existence of optimal solutions for continuous-time OCPs, one will typically invoke further differentiability properties of the problem data and feasibility assumptions. For the sake of brevity, we skip the details and refer to classical texts such as \citep{Lee67,Berkovitz74} or more recent treatments \citep{Hartl95,Vinter10,Gerdts11a}.

Comparing the discrete-time and the continuous-time variants---OCP \eqref{eq:ocp_dt}, respectively, OCP \eqref{eq:ocp_ct}---shows that as OCP \eqref{eq:ocp_dt} is defined in normed metric spaces it admits the interpretation as an infinite-dimensional system. This implies that a large variety of discrete-time results to be touched upon later can be---or already have been---formulated in this general setting. In contrast, it is more difficult to formalize infinite-dimensional continuous-time optimal control problems in a unified way. Hence we have limited \eqref{eq:ocp_ct} to finite-dimensional problems. Nevertheless, we will comment on infinite-dimensional continuous-time problems at several places throughout this chapter.

\subsection{Optimality Conditions}
Dual variables are key to the understanding of any optimization problem. Hence, we now formulate the first-order Necessary Conditions of Optimality (NCOs) for OCP \eqref{eq:ocp_dt}, respectively, OCP \eqref{eq:ocp_ct}.

The NCOs for OCP \eqref{eq:ocp_dt} are built using the point-wise Lagrangian function
\begin{subequations} \label{eq:L_dt}
	\begin{align}
	\mcl{L}(0) &\doteq \lambda_0^\top(x(0)-x_0), \\
	\mcl{L}(t) &\doteq \ell(x(t),u(t)) +\lambda(t+1)^\top(f(x(t), u(t)) -x(t+1)) \nonumber \\
	&\phantom{\doteq} +\mu(t)^\top g(x(t), u(t)), \qquad t \in   \mathbb{N}_{[0,T-1]},\\
	\mcl{L}(T) &\doteq \varphi(x(T)) +\mu(T)^\top \psi(x(T)),
	\end{align}
\end{subequations}
which directly leads to the overall Lagrangian 
\[
{L}(x(\cdot),u(\cdot),\lambda(\cdot),\mu(\cdot)) \doteq \sum_{t=0}^T \mcl{L}(t),
\]
with $\mcl{L}(t)$ from \eqref{eq:L_dt}. Here, $x$ and $u$ are called the primal variables while $\lambda$ and $\mu$ are called the dual variables or adjoints/costates (in case of $\lambda$).

Let the objective functional in OCP \eqref{eq:ocp_dt} be written as
\[
J_T(u(\cdot), x_0) = \sum_{t=0}^{T-1} \ell(x(t),u(t))  + \varphi(x(T)) 
\]
and the one in OCP \eqref{eq:ocp_ct} as
\[
J_T(u(\cdot), x_0) =  \int_{0}^{T} \ell(x(t),u(t))\mathrm{d}t  + \varphi(x(T)). 
\]
Then, for finite horizons $T < \infty$, optimality of primal solutions $(x(t),u(t))$  is characterized by
\[
J_T(u^\star(\cdot), x_0) \leq J_T(u(\cdot), x_0) \quad \text{ for all admissible } u(\cdot)
\]
A dual solution is called optimal if it satisfies the NCOs defined below together with the primal optimal solution. Optimal primal and dual solutions, provided they exist, are denoted by superscript $\cdot^\star$, i.e.\ we write  $\chi^\star(\cdot;x_0), \,\chi \in\{x, u, \lambda, \mu\}$. The argument $x_0$ will be dropped whenever the considered initial condition can be inferred from the context. 
The first-order NCOs are given by $\nabla {L}=0$, which entails 
\begin{subequations} \label{eq:NCO_dt}
\begin{align}
\hspace*{-5mm}x^\star(t+1) &= f(x^\star(t), u^\star(t)),\hspace*{2.03cm} x^\star(0) = x_0, \\
 \lambda^\star(t) &=  f_x^\top\lambda^\star(t+1) +\ell_x+ g_x^\top\mu^\star(t), ~
 \lambda^\star(T) =  \left(\varphi_x + \psi_x^\top \mu^\star\right)\big|_{t=T}, \label{eq:NCO_dt_trans}\\
 0 &= f_u^\top\lambda^\star(t)+\ell_u +g_u^\top\mu^\star(t).
\end{align}
\end{subequations}
Observe that the NCOs as listed above leave out parts of standard KKT conditions such as primal feasibility, complementary slackness, or higher-order conditions. However, for the purposes of the subsequent exposition, the conditions above suffice. 

In correspondence to the discrete-time OCP we consider the steady-state optimization problem 
\begin{subequations}
	\label{eq:sop_dt}
	\begin{align}
\min_{\bar u, \bar x} \ \ &  \ell(\bar x,\bar u)   \label{eq:sop_dt_obj} \\
	\text{subject to} \nonumber \\  
	 \bar x &= f(\bar x,\bar u),  \label{eq:sop_dt_dyn}\\
	 0&\geq g(\bar x,\bar u) \label{eq:sop_dt_pc}
	\end{align}
\end{subequations}
and its Lagrangian
$
\bar L = \ell(\bar x,\bar u) + \bar{\lambda}^\top \left(f(\bar x, \bar u) -\bar x\right)+ \bar{\mu}^\top g(\bar x, \bar u)$. Superscript $\bar{\cdot}$ denotes variables at steady state and hence optimal steady-state variables are indicated by $\bar{\cdot}^\star$.\footnote{A feasible primal pair $(\bar x, \bar u)$ satisfying \eqref{eq:sop_dt_dyn} in discrete time or \eqref{eq:sop_ct_dyn} in continuous time is called a \textit{controlled equilibrium}.\label{fn:coneq}}
The optimality  conditions of this problem entail
\begin{subequations} \label{eq:NCO_s_dt}
\begin{align}
\bar x^\star &= f(\bar x^\star, \bar u^\star),\\
 \bar \lambda^\star &=  f_x^\top\bar\lambda^\star +\ell_x+ g_x^\top\bar\mu^\star,\\ 
 0 &= f_u^\top\bar \lambda^\star+\ell_u +g_u^\top\bar\mu^\star.
\end{align}
\end{subequations}
Notice that the optimality conditions \eqref{eq:NCO_s_dt} of the stationary problem \eqref{eq:sop_dt} coincide with the full NCOs \eqref{eq:NCO_dt} of OCP \eqref{eq:ocp_dt} if only constant solutions are considered.\vspace*{2mm}

The continuous-time counterpart to \eqref{eq:NCO_dt} can be obtained via a classical direct adjoining approach, i.e.\ one considers the Hamiltonian
\begin{equation}\label{eq:H}
H(x,u,\lambda,\mu) \doteq \lambda^0\ell(x,u) +\lambda^\top f(x,u) + \mu^\top g(x,u).
\end{equation}
Excluding abnormal problems we set $\lambda^0 \equiv 1$ and then the Pontryagin Minimum Principle  reads
\begin{subequations} \label{eq:NCO_ct}
\begin{align}
\dot x^\star &= \phantom{-}f(x^\star, u^\star),  \hspace*{1.48cm} x^\star(0) = x_0, \\
\dot \lambda^\star &=  -f_x^\top\lambda^\star -\ell_x- g_x^\top\mu^\star,~
 \lambda^\star(T) =  \left(\varphi_x +\psi_x^\top\mu^\star\right)\big|_{t=T},  \label{eq:NCO_ct_trans}\\
 0 &= \phantom{-}f_u^\top\lambda^\star+\ell_u +g_u^\top\mu^\star,\\
 &H(x^\star,u^\star,\lambda^\star,\mu^\star) \leq \min_{u\in \mbb{U}(x^\star)}H(x^\star,u,\lambda^\star,\mu^\star).
\end{align}
\end{subequations}
Similar to the discrete-time case, we have left out some technical conditions which do not affect the further developments. Specifically, in the time-invariant setting considered here the Hamiltonian is constant along optimal solutions $H(x^\star,u^\star,\lambda^\star,\mu^\star) = const$ and the adjoints satisfy $(\lambda_0, \lambda^\star(t)) \not= 0$ for all $t\in [0, T]$. Moreover, in analyzing the continuous-time OCPs state path constraints \eqref{eq:ocp_ct_pc} can lead to difficulties. Hence, one usually either restricts the setting to mixed-input-state constraints combined with a linear independence condition, or one limits the degree of pure state path constraints, or one considers optimality conditions formulated in terms of bounded variation. We refer, e.g., to \cite{Hartl95} for comments and insights. Here, however, we will be able to avoid those technicalities by assuming for OCP \eqref{eq:ocp_dt}, respectively, OCP \eqref{eq:ocp_ct} that the turnpike properties to be analyzed occur in the interior of the constraint set $\mbb{Z}$.

The steady-state optimization problem corresponding to OCP \eqref{eq:ocp_ct} reads
\begin{subequations}
	\label{eq:sop_ct}
	\begin{align}
\min_{\bar u, \bar x} \ \ &  \ell(\bar x,\bar u)   \label{eq:sop_ct_obj} \\
	\text{subject to} \nonumber \\  
	 0 &= f(\bar x,\bar u),  \label{eq:sop_ct_dyn}\\
	 0&\geq g(\bar x,\bar u). \label{eq:sop_ct_pc}
	\end{align}
\end{subequations}
The reader will face no difficulties in verifying that also in continuous time the optimality conditions  of the stationary problem \eqref{eq:sop_ct} correspond to the full NCOs \eqref{eq:NCO_ct} of OCP \eqref{eq:ocp_ct} if constant solutions are considered.

\section{Definition and Taxonomy of Turnpike Properties}\label{sec:Taxonomy}

In many applications it is of interest to consider the solutions of OCPs for varying problem data such as parameters of the underlying dynamics---\eqref{eq:ocp_dt_dyn} respectively \eqref{eq:ocp_ct_dyn}---, different initial conditions $x_0$ of said model, or varying  length of the optimization horizon $T$ in \eqref{eq:ocp_dt}, respectively, in \eqref{eq:ocp_ct}. Classical sensitivity theory \citep{Fiacco76} provides a handle for locally analyzing the former, i.e.\ for small parametric variations. In contrast, turnpike theory aims at a \textit{non-local} analysis of the structure of  solutions for varying initial conditions $x_0$ and horizons $T$ over a large range of values.  Hence, we consider OCP \eqref{eq:ocp_dt} / \eqref{eq:ocp_ct} for a set of initial conditions $\mbb{X}_0\subseteq\mbb{X}\subseteq\mcl{X}$ and horizons $T\in \mbb{R}^+/\mbb{N} $. 
Notice that, depending on the considered input constraints and the properties of \eqref{eq:ocp_dt_dyn} / \eqref{eq:ocp_ct_dyn}, the set $\mbb{X}_0$ may be open or closed and bounded or unbounded. However, in most relevant cases it will be uncountable. As we will see later, the size of the maximal set $\mbb{X}_0$ depends on the reachability properties of \eqref{eq:ocp_dt_dyn} / \eqref{eq:ocp_ct_dyn}.

Recall the notational convention for optimal primal and dual solutions $\chi^\star(\cdot;x_0)$, $\chi \in\{x, u, \lambda, \mu\}$. Let $\bar \chi$ be some controlled equilibrium satisfying \eqref{eq:sop_dt_dyn} / \eqref{eq:sop_ct_dyn}  (see Footnote \ref{fn:coneq} for an informal definition), and let $\mcl{B}_\varepsilon(\bar\chi)$ be a ball of radius $\varepsilon > 0$ centred at $\bar \chi$.
Turnpike theory tries to establish global properties of the following type:
\begin{equation}
\forall x_0 \in \mbb{X}_0, \forall T >0, \forall \varepsilon >0,  ~|~\left[ \text{time }\chi^\star(\cdot;x_0)\not\in \mcl{B}_\varepsilon(\bar \chi)\right] \leq \nu(\varepsilon) < \infty
\label{eq:tpinformal}\end{equation}
where $\nu: [0,\infty) \to \mbb{R}\cap\infty$ is independent of the actual horizon length $T$ and of the initial condition $x_0$. Naturally, the conceptual property above calls for precisions and further comments: With respect to which variables $\chi \in\{x, u, \lambda, \mu\}$ shall the closeness be established? How to measure/express the time spent outside of $\mcl{B}_\varepsilon(\bar \chi)$? Besides these technical aspects it is imperative to notice two issues: (a) The conceptual statement above becomes trivial, if the considered horizons $T$ are bounded from above by a finite number. (b) The property is to be established for all $x_0 \in \mbb{X}_0$. Thus, given $\mbb{X}_0$ is more complex than a simple $\varepsilon$-ball, one aims at establishing a non-local (similarity) property of optimal solutions for varying initial conditions and varying horizon length. 

\subsection{Motivating Example} \label{sec:Fish}
To illustrate the turnpike phenomenon, we consider a simple example dating back to \cite{Cliff73}. 
The OCP at hand reads
\begin{subequations} \label{eq:Fish}
\begin{align} 
\min_{ u(\cdot) \in \mcl{L}([0,T], \mcl{U}) } 
&\int_{0}^T \big[ax(t) + bu(t) -cx(t)u(t)\big]\,\mathrm{d} t, \label{eq:objFish}\\
\text{subject to}& \nonumber \\
\dot x &= x(x_s-x-u), \quad x(0) = x_0 > 0 \label{eq:dynFish}\\
\forall t \in [0,T|: u(t) &\in [0, u_{max}] \\
\forall t \in [0,T|: x(t) &\geq x_{min}
\end{align}
\end{subequations}
The state $x$ is the fish density in a certain habitat, the control $u$ is the fishing rate, and
the parameter $x_s$ describes the highest sustainable fish density. In slightly modified form---i.e., with an additional terminal constraint and the state constraint $x(\tau) > 0$---this OCP is analyzed by
\cite{Cliff73} and in \cite[Chap. 3.3]{CaHL91}. 

The parameter values are $a=1, b=c=2, x_s = u_{max} = 5, x_{min} = 0.1$. 
In case the state constraint  $x\geq x_{min}$ is active on some interval $[\tau_{act}, \tau_{deact}]$, it is easily verified that 
$u(t) = x_s - x_{min} >0, \quad \forall t \in [\tau_{act}, \tau_{deact}]$.
However, for sufficiently small but positive $x_{min}$,  if $x_0 = x_{min}$, the optimal solutions depart from the state constraint immediately, 

The considered OCP has exactly one singular arc of order two along which 
\begin{equation}\label{eq:Ex_sing_arc}
\bar x = \frac 12 x_s +\frac{b-a}{2c} > x_{min}, \quad \bar u = x_s -\bar x, \quad \bar \lambda = c-\dfrac{b}{\bar x}, \quad \bar\mu = 0,
\end{equation}
i.e.\ primal and dual variables are at steady state.
We refer to \citep{Bryson69a,Zelikin94} for details on the analysis of singular arcs and to \cite[Chap. 3.3]{CaHL91} for details on the derivation of the singular arc in this example.
As long as the state constraint does not become active, the optimal input signal for this OCP takes only three values
$u^\star(t; x_0) \in \{0,\bar u, u_{max}\}$ for all $t\in [0,T]$ and all $x_0 \geq x_{min}$.
 
\begin{figure}[t]
 \includegraphics[scale=0.375]{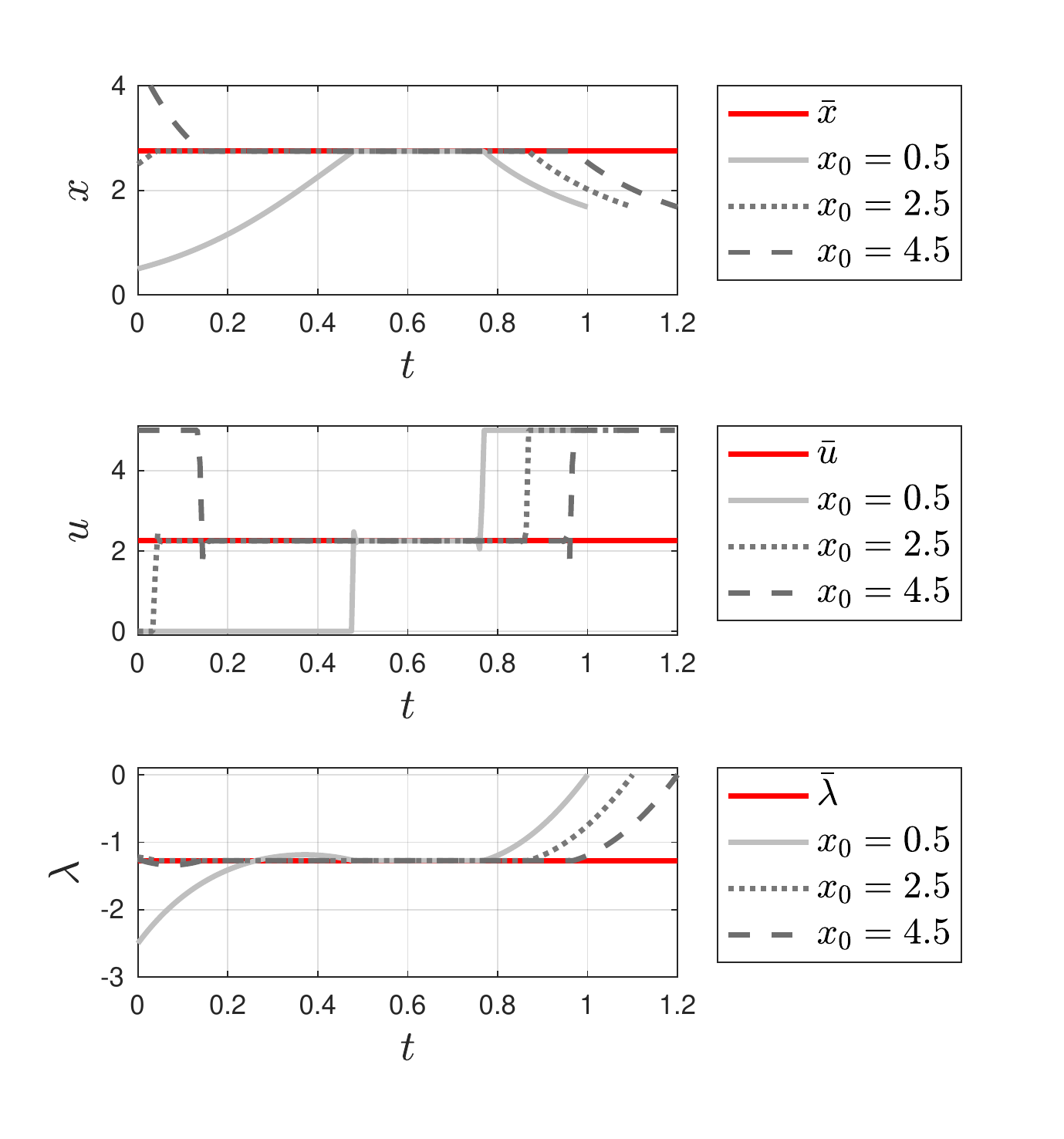} 
 \caption{Fish harvest example with bilinear objective with pairs of $x_0$ and $T$: $x_0 = 0.5, T = 1.0$;  $x_0 = 2.5, T = 1.1$; and  $x_0 = 4.5, T = 1.2$. Depicted are optimal state, input, and adjoint trajectories. }\label{fig:Ex_singular}
 \end{figure} 
 
 Figure \ref{fig:Ex_singular} shows numerical results for initial conditions $x_0 \in\mbb{X}_0 = \{0.5,~2.5,$ $4.5\}$ and horizons $T \in \{1.0,~1.1,~1.2\}$. As one can see, the optimal solutions all comprise three distinct arcs: the turnpike entry arc, the turnpike arc, and the leaving arc. Put differently, all solutions enter the steady state at $(\bar x, \bar u)$ exactly.\footnote{We remark that for turnpike analysis it is in general \textit{not crucial} that the optimal solutions enter the turnpike at $(\bar x, \bar u)$ \textit{exactly}. Indeed, as we comment subsequently, in many cases the solutions will only get arbitrarily close, see also Remark \ref{rem:exactTP} and the specific turnpike definitions below.} Moreover, while the primal variables enter  $(\bar x, \bar u)$, the adjoint $\lambda$ also reaches its steady value $\bar \lambda$. 
Notice that, despite different initial conditions and different horizons, the optimal solutions are qualitatively similar, i.e.\ they have the middle arc along which the optimal solutions are stationary in common.  As demanded by \eqref{eq:tpinformal}, the time the optimal solutions spent outside of the steady state triplet $(\bar x, \bar u, \bar \lambda)$ is bounded independent of the horizon length.  

\subsection{Turnpike Definitions}
In order to be more precise, we turn towards a discrete-time definition of the turnpike phenomenon which was introduced by \cite{Grue13}, see also \cite[first paragraph of Chapter 8]{Zasl06} for an earlier though different definition. To this end, consider the placeholder variable $\xi \in \{x, (x,u), (x,u,\lambda, \mu)\}$, which refers to time-dependent functions of the state, input, and dual variables and its steady-state counterpart $\bar \xi$.

\begin{definition}[Discrete-time (cardinality) turnpike]\label{def:CardTP}
OCP~\eqref{eq:ocp_dt} is said to have a turnpike at $\bar \xi, ~\xi \in \{x, (x,u), (x,u,\lambda, \mu)\}$ if there exists a function $\nu_\xi:[0,\infty)\to [0,\infty]$ such that, for all $x_0 \in \mbb{X}_0$ and all $T\in \mbb{N}$,
	\begin{equation}\label{eq:TP_dt}
	\#\mathbb{Q}_{\xi,T}(\varepsilon)  \leq\nu_\xi(\varepsilon)< \infty \quad \forall\: \varepsilon >0,
	\end{equation}
 where 
	\begin{equation}\label{eq:SetQ}
	\mathbb{Q}_{\xi,T}(\varepsilon)\doteq\{t\in\{0,\dots,T-1\}\hspace*{1mm}|\hspace{1mm} \|\xi^\star(t;x_0)-\bar{\xi}\| > \varepsilon \}
	\end{equation}
	and $\#\mbb{Q}_{\xi,T}(\varepsilon)$ is the cardinality of $\mbb{Q}_{\xi,T}(\varepsilon)$.
	\begin{itemize}
	\item If $\xi = x$, then OCP~\eqref{eq:ocp_dt} is said to have a \textit{state turnpike} at $\bar x$.
	\item If $\xi = (x, u)$, then OCP~\eqref{eq:ocp_dt} is said to have an \textit{input-state turnpike} at $(\bar x, \bar u)$.
	\item If $\xi = (x, u, \lambda, \mu)$, then OCP~\eqref{eq:ocp_dt} is said to have a \textit{primal-dual turnpike} at $(\bar x, \bar u, \bar 
	\lambda, \bar \mu)$.
	\end{itemize}
	Moreover, if 
	\[
	\lim_{\varepsilon \shortarrow{7} 0}\nu_\xi(\varepsilon) < \infty
	\]
	then the turnpike property is said to be \textit{exact}. 
\end{definition}
Usually, one refers to the steady-state variables $\{\bar x, (\bar x,\bar u), (\bar x, \bar u,\bar \lambda, \bar \mu)\}$ simply as the \textit{turnpike}. Historically, many results have addressed the case of state turnpikes, while the input-state and the primal-dual case have received attention more recently. One important aspect in these developments has been the research on converse turnpike results, i.e.\ analyzing the implications if an OCP is assumed to exhibit a specific turnpike property. We will come back to this point later in Section \ref{sec:DI_TP}.

The continuous-time counterpart to Definition \ref{def:CardTP} was proposed in \citep{epfl:faulwasser14e,epfl:faulwasser15g}. However, the conceptual idea has already been used in earlier results \cite[Thm. 4.2]{CaHL91}. 
Using again the placeholder variable $\xi$, we define the set 
\begin{equation} \label{eq:SetTheta}
\Theta_{\xi,T}(\varepsilon) \doteq \left\{t \in [0,T]\,|\, \left\| \xi^\star(t; x_0) - \bar \xi\right\| > \varepsilon\right\}.
\end{equation}

\begin{definition}[Continuous-time (measure) turnpike]\label{def:MeasTP}
OCP~\eqref{eq:ocp_ct} is said to have a turnpike at $\bar \xi, \xi\in \{x, (x,u), (x,u,\lambda, \mu)\}$ if there exists a function $\nu_\xi:[0,\infty)\to [0,\infty]$ such that, for all $x_0 \in \mbb{X}_0$ and all $T \in \mbb{R}^+$, we have
\begin{equation}\label{eq:TP_ct}
\mu\left[\Theta_{\xi,T}(\varepsilon)\right]  \leq\nu_\xi(\varepsilon)< \infty \quad \forall\: \varepsilon >0,
\end{equation}
where $\mu[\cdot]$ is the Lebesgue measure on $\mbb{R}$.
\begin{itemize}
	\item If $\xi = x$, then OCP~\eqref{eq:ocp_ct} is said to have a \textit{state turnpike} at $\bar x$.
	\item If $\xi = (x, u)$, then OCP~\eqref{eq:ocp_ct} is said to have an \textit{input-state turnpike} at $(\bar x, \bar u)$.
	\item If $\xi = (x, u, \lambda, \mu)$, then OCP~\eqref{eq:ocp_ct} is said to have a \textit{primal-dual turnpike} at $(\bar x, \bar u, \bar 
	\lambda, \bar \mu)$.
	\end{itemize}
Moreover, if 
	\[
	\lim_{\varepsilon \shortarrow{7} 0}\nu_\xi(\varepsilon) < \infty
	\]
	then the turnpike property is said to be \textit{exact}.
\end{definition}

The two previous definitions (in their non-exact variants) state the most general turnpike properties we will encounter in this chapter. The reason being that, in principle, the measure/cardinality based properties allows for cases where optimal solutions would temporarily depart and return from the turnpike neighborhood. For the case of exact state turnpikes it is easy to verify that this cannot be the case \cite[Lem.\ 1]{kit:faulwasser17a}.\footnote{Occasionally, exact turnpike are also denoted as \textit{finite} turnpikes \citep{Gugat20a}. Here, however, we stick to the wording coined in \cite[Chap. 2]{CaHL91}.}
For the non-exact case, which was called approximate in \citep{epfl:faulwasser15g},  similar statements can be made using dissipativity conditions \citep{FGHS20}. 

To continue our exploration of different turnpike concepts, we turn toward a stronger notion where instead of measures and cardinalities the deviation of optimal solutions from their turnpike values is bounded in time by an exponential function. This kind of turnpike property was first observed in the economics literature, see, e.g., \cite[Theorem (4.5)]{Bewl82} or \cite[Theorem 10.1]{McKe86}. Here we give a definition that appears in similar form, e.g., in \citep{PorZ13} or \citep{DGSW14}.

\begin{definition}[Exponential turnpike] \label{def:ExpTP}
The discrete-time OCP~\eqref{eq:ocp_dt} is said to have an exponential turnpike at $\bar \xi, ~\xi\in \{x, (x,u), (x,u,\lambda, \mu)\}$ if there exist $C >0, \rho\in [0,1)$ such that,
for all $x_0 \in \mbb{X}_0$ and all $T\in \mbb{N}$,
	\begin{equation}\label{eq:expTP_dt}
	\|\xi^\star(t;x_0)-\bar{\xi}\| \leq C\left(\rho^{t}+ \rho^{T-t} \right).
	\end{equation}	
Similarly, the continuous-time OCP~\eqref{eq:ocp_ct} is said to have an exponential turnpike at $\bar \xi, \xi\in \{x, (x,u), (x,u,\lambda, \mu)\}$ if there exist $C >0, \gamma >0$ such that,
for all $x_0 \in \mbb{X}_0$ and all $T\in \mbb{R}$, if \eqref{eq:expTP_dt} holds with $\rho \doteq e^{-\gamma}$.
\end{definition}

One can show that the exponential turnpike property implies the cardinality/measure variant as one computes from \eqref{eq:expTP_dt} that 
\[
\mu\left[\Theta_{\xi,T}(\varepsilon)\right] = \max\left\{0, 2t_1(\varepsilon, T)       \right\},\quad t_1(\varepsilon, T) = -\dfrac{1}{\gamma} \ln\left(\frac{\varepsilon}{2C} + \sqrt{\frac{\varepsilon^2}{4C^2}-e^{-\gamma T}}\right).
\]
Here $t_1(\varepsilon, T)$ can be interpreted as the worst-case time the optimal solutions enter the $\varepsilon$-ball around the turnpike. This derivation uses that the symmetry of the right hand side in \eqref{eq:expTP_dt} implies $T-t_1(\varepsilon, T) = t_2(\varepsilon, T)$, where $t_2(\varepsilon, T)$ is the worst-case time the optimal solutions leave the $\varepsilon$-ball.
Hence, one has that
\[
\mu\left[\Theta_{\xi,T}(\varepsilon)\right] \leq \max \left\{0, -\frac{2}{\gamma} \ln\frac{\varepsilon}{2C}\right\}.
\] 
Similarly, in the discrete-time case  \eqref{eq:expTP_dt} implies
\[
\#\mbb{Q}_{\xi,T}(\varepsilon) \leq \min\left\{k\in\mbb{N}\,\big|\,-\frac{2}{\gamma} \log_\rho\frac{\varepsilon}{2C}\leq k\right\}.
\]
Summing up, the definitions of the cardinality/measure turnpike and the exponential variants provide the blueprint for a number of specific variants such as, e.g., exponential input-state turnpikes etc.
\begin{figure}[t]
\begin{center}
\scalebox{0.75}{
\begin{tikzpicture}	
	  \node (Title) at (-0.75,.75)   {Characterization by Asymptotics };
		  \draw [draw=black] (1.75,-1.95) rectangle (2.55,0.3);
		  \node (Aprim) at (2.15,-.8){ \rotatebox{90}{\small exact turnpike }};
		  \draw [draw=black] (-3.45,-1.05) rectangle (2.75,0.4);
		  \node (Air) at (-.5,-.35) {exponential turnpike  };		
	  	\draw [draw=black] (-4.45,-2.25) rectangle (2.95,0.5);
		  \node (Airr) at (-.5,-1.65) {measure/cardinality turnpike}; 
		  \node (shift) at (0, -4) {~};
		  	\end{tikzpicture}
		  	\quad
		\begin{tikzpicture}
			  \node (Title) at (0,.75)   {Characterization by Argument };
			 \draw \boundellipse{-1.35,-0.75}{2.25}{0.75};
			 \node (sTP) at (-2.15,-.75) {state turnpike  };
			\draw \boundellipse{ 1.35,-0.75}{2.25}{0.75};
	  \node (iTP) at (2.15,-.75) {input turnpike  }; 		
		  \node (isTP) at (0.,-.75) {\small{input-state}  };	
		  \draw  \boundellipse{ 0,-2.95}{.75}{2.05};
		   \node (pdTP) at (-.0,-2.85) {\rotatebox{90}{dual turnpike}};	 
		   \draw[-] (0,-1.1) -- (1.25,-1.95);
		  \node (pdTP) at (2.25,-1.95) {primal-dual};	 
		    \node (pdTP) at (2.25,-2.35) {turnpike};	 
		\end{tikzpicture}
		}
\end{center}
	\caption{Taxonomy of turnpike properties with respect to the asymptotics of the optimal solutions (left) and with respect to considered arguments (right).
	\label{fig:taxo} }
\end{figure}
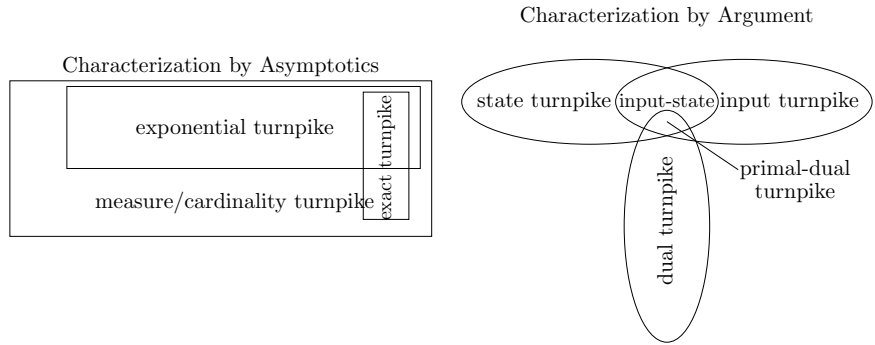
That is, one can create taxonomy of turnpike properties with respect to two different dimensions:
\begin{itemize}
\item characterization by the asymptotics  of \textit{how} the closeness of optimal solutions to the turnpike is measured---i.e.\ the difference of Definition \ref{def:CardTP}/\ref{def:MeasTP} and Definition \ref{def:ExpTP}---, and
\item characterization by the arguments/variables of \textit{what} is close to its turnpike value---i.e., the difference between subvariants in Definition \ref{def:CardTP}/\ref{def:MeasTP}.
\end{itemize}
This taxonomy is sketched in Figure \ref{fig:taxo}
and we notice that it
is not entirely complete. For example  \cite{Gugat19a} investigate a stronger-than-measure-turnpike weaker-than-exponential-turnpike result for OCPs comprising infinite-dimensional hyperbolic systems and strict convexity assumptions for the objective. Specifically, it is shown that the optimal control  profiles are close to their stationary values in a time-averaged sense. In other words, instead of an exponential bound---see \eqref{eq:expTP_dt} for the discrete-time version---a weaker integrable bound is established. Theorem 1 and Theorem 2 in \citep{Gugat19a}  establish input turnpike properties. The latter under integrality constraints for the controls. 
 For problems arising in OCPs comprising finite-dimensional mechanical systems,  \citep{kit:faulwasser19b_2,tudo:faulwasser20e} recently proposed the consideration of hyperbolic bounds. \vspace*{2mm}

During the first 50-60s years after the foundational quote by \cite{Dorfman58},---and quite surprisingly given the numerous works and the impact of turnpike results---precise definitions of turnpike properties have been, to the best of our knowledge, not been common in the literature. Actually, the book of \cite{Zasl06} appears to be the first to present an explicit turnpike definition. This implies that the vast majority of classical turnpike results does not adhere to any kind of taxonomy and define the property under consideration in an ad-hoc manner. Rather most of those results are denoted simply as \textit{turnpike results}. 

Putting the taxonomy proposed above in use, we note that originally state turnpike properties have been in the focus of many works in economics \citep{Mckenzie76,CaHL91}. 
 The need to extend to input-state turnpikes arose in receding-horizon and predictive control \citep{epfl:faulwasser15g} and in the investigation of converse turnpike results \citep{GruM16,FKJB17}. We remark that especially converse turnpike analysis relies on precise definitions of the phenomenon.
 In parallel to the developments driven by receding-horizon and predictive control, the move towards primal-dual (exponential) turnpikes was done by \cite{PorZ13,TreZ15}. We will comment in Section \ref{sec:DI_TP} on the differences in the technical analysis of the varying approaches. Moreover, in Section \ref{sec:infHor} we will comment on the extension towards infinite horizons.

\section{Generating Mechanisms} \label{sec:Results}
Before entering the technical analysis we informally illustrate the mechanisms that generate turnpike phenomena. To this end, we consider the three distinct arcs of turnpike solutions---the entry arc, the turnpike arc, and the leaving arc---separately. Notice that the characteristic leaving arc is not required by the Definitions \ref{def:CardTP}--\ref{def:ExpTP}. 

For the sake of illustration, we revisit the fish harvest problem from Section \ref{sec:Fish}.
\begin{figure}[htb]
 \includegraphics[scale=0.375]{figures/t_sing.pdf}   \includegraphics[scale=0.375]{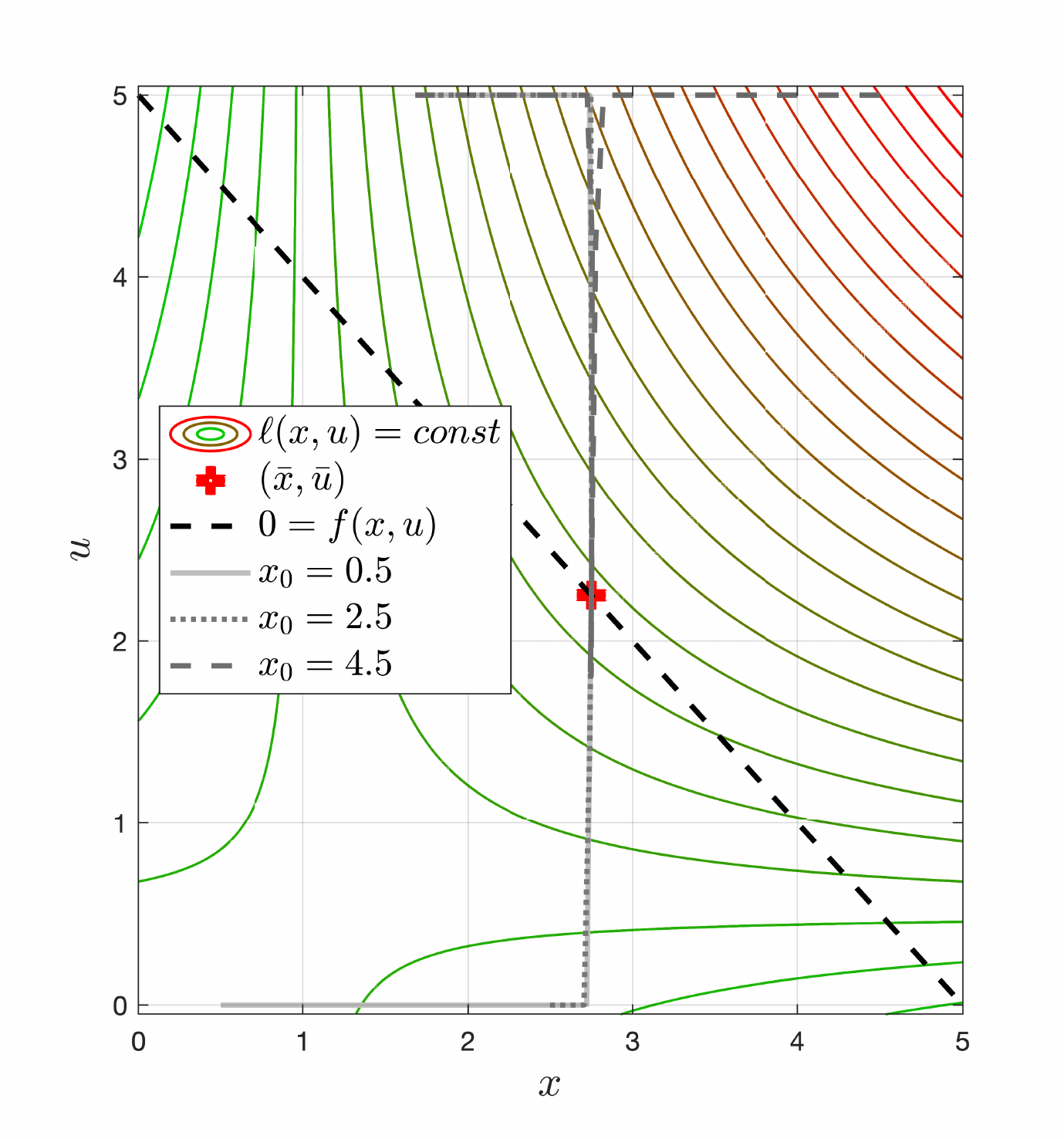}
 \caption{Fish harvest example with bilinear objective with pairs of $x_0$ and $T$: $x_0 = 0.5, T = 1.0$;  $x_0 = 2.5, T = 1.1$; and  $x_0 = 4.5, T = 1.2$.  Left: optimal state, input and adjoint trajectories. Right: $x-u$ plane with optimal trajectories and contour lines of the Lagrange function $\ell(x,u)$ (coloured), and the set of admissible steady states (dashed black). }\label{fig:Ex_singularII}
 \end{figure}  
 Figure \ref{fig:Ex_singularII} illustrates the geometry of OCP \eqref{eq:Fish} in the $x-u$ plane (right) while also recalling the optimal trajectories (left). The dashed black line in the right hand side figure depicts the set of steady states given by $0 = x(x_s-x-u),$ and the coloured lines are the contour lines of $\ell(x,u)$. The red cross marks the optimal steady state  at which the set of steady states is tangent to the contour lines of the objective. This point coincides with the primal variables of the singular arc \eqref{eq:Ex_sing_arc}.

Figure \ref{fig:Ex_singularII} provides an intuitive explanation why the problem at hand exhibits the turnpike phenomenon: It is easily verified that the singular arc \eqref{eq:Ex_sing_arc} corresponds to the unique KKT point of the steady-state optimization problem \eqref{eq:sop_ct} for OCP  \eqref{eq:Fish}. Moreover, one can directly infer that above the steady-state set ($=$ dashed black line), the dynamics of \eqref{eq:dynFish} points in the direction of lower $x$ values, while below the steady-state set it points in the direction of increasing $x$ values.

In other words, the fact that the optimal solutions approach the turnpike clearly requires (at least asymptotic) controllability properties of the dynamics. Moreover, reaching the turnpike implicitly expresses that there does not exist any other set or state than the turnpike $(\bar x, \bar u)$ where the system can stay infinitely long {and} which has better time-averaged performance than $\ell(\bar x, \bar u)$.\footnote{This observation indicates that the turnpike phenomenon is closely linked to a property which was coined at optimal operation at steady state, which in turn is closely linked to overtaking optimality \citep{CaHL91}. This property requires that in a time-averaged sense a dynamic system is optimally operated at a specific steady state. For details and the link to dissipativity of OCPs see \citep{Mueller14a} for the discrete-time case and \citep{FKJB17} for continuous-time results. }

Finally, comparing the left and right hand side plots in Figure \ref{fig:Ex_singularII}, one may observe that the characteristic leaving arc steers the system to the transient domain above the steady-state set, where for the final part of the time horizon the optimal solutions are at a domain of lower values of $\ell$. Note that, however, the system can only stay a limited amount of time in this domain due to the system dynamics and the constraints.
\vspace*{2mm}

\begin{figure}[htb]
 \includegraphics[scale=0.375]{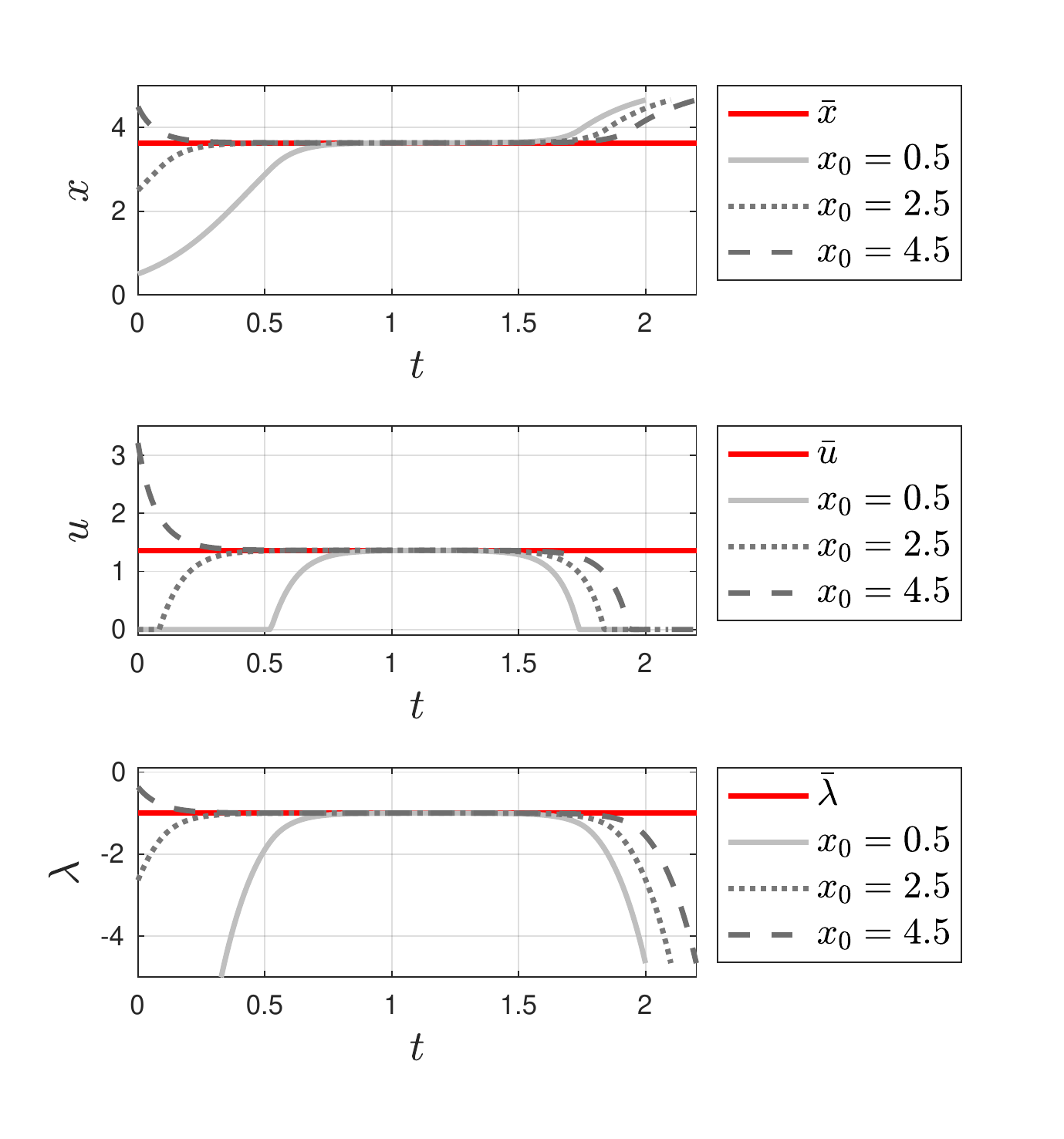}   \includegraphics[scale=0.375]{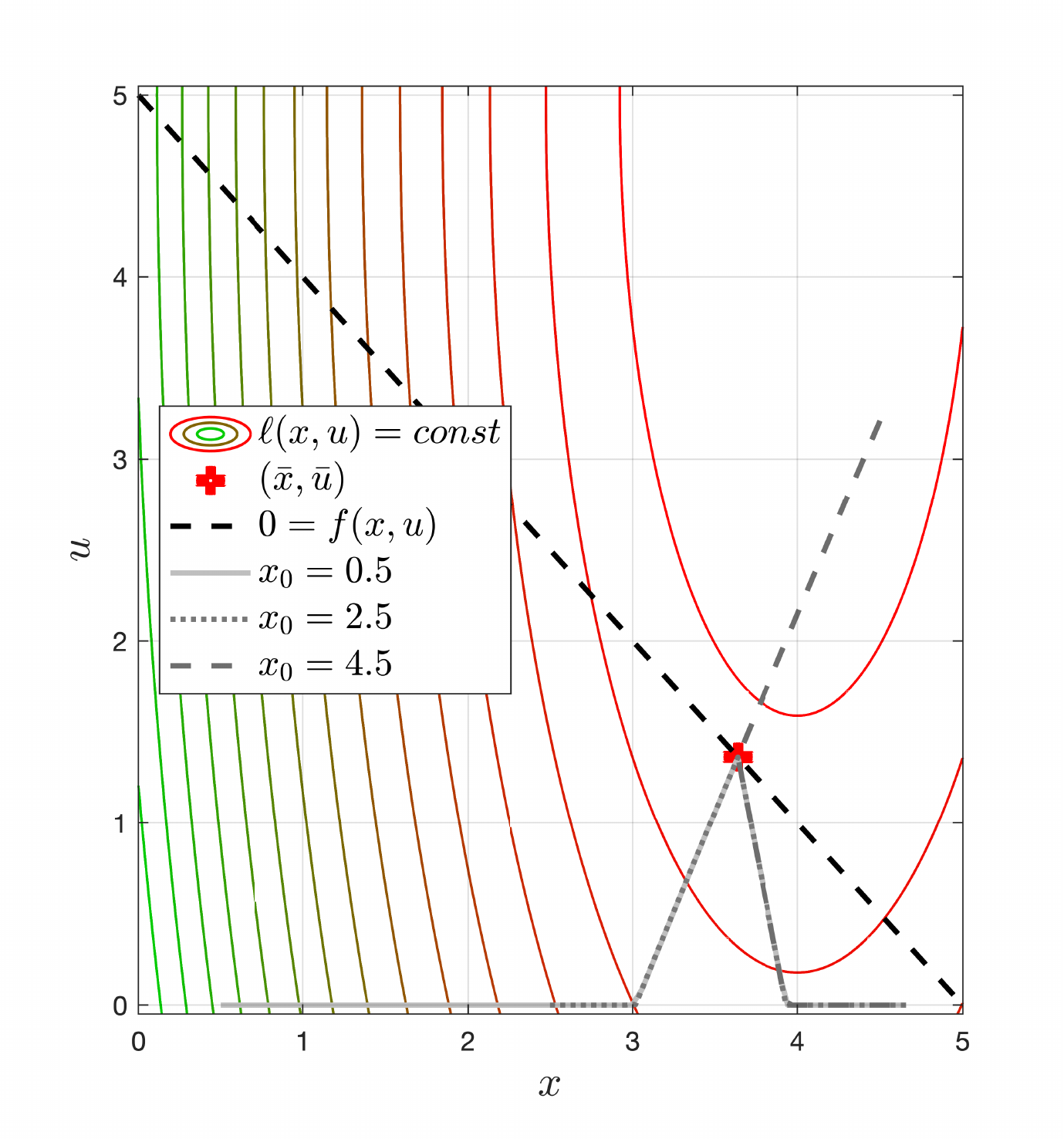}
 \caption{Fish harvest example with quadratic objective and Mayer term \eqref{eq:Fish_quadI} with pairs of $x_0$ and $T$: $x_0 = 0.5, T = 2.0$;  $x_0 = 2.5, T = 2.1$; and  $x_0 = 4.5, T = 2.2$. Left: optimal state, input and adjoint trajectories. Right: $x-u$ plane with optimal trajectories and contour lines of the Lagrange function $\ell(x,u)$ (coloured), and the set of admissible steady states (dashed black). }\label{fig:Ex_quadI}
 \end{figure}

At first glance, one may conjecture that the characteristic turnpike phenomenon is induced by the fact that for sufficiently long horizons OCP \eqref{eq:Fish} exhibits the unique and stationary singular arc \eqref{eq:Ex_sing_arc}. 
However, as we show next, this is not the case.
To this end, we modify OCP \eqref{eq:Fish} as follows
\begin{equation} \label{eq:Fish_quadI}
\ell(x,u) = \dfrac{1}{2}\left( (x - x_c)^2q +  (u - u_c)^2\right), \quad \varphi(x) = -\dfrac{1}{2} x^2,\quad  \psi(x) = 0
\end{equation}
with $x_c = 4, u_c = 5$ and $q = 10$. Due to the quadratic input penalization, OCP \eqref{eq:Fish} with \eqref{eq:Fish_quadI} is not singular.  The left hand side of Figure \ref{fig:Ex_quadI} shows its numerical solutions. As one can see, despite the changes to the objective, the optimal solutions also show the turnpike phenomenon. The right hand side of Figure \ref{fig:Ex_quadI}  illustrates again the geometry in the $x-u$ plane.
Recall that without a Mayer term (and with bilinear Lagrange term $\ell$) the optimal solutions of OCP \eqref{eq:Fish} exhibited a leaving arc approaching the domain with lower values of $\ell$, cf.\ Figure \ref{fig:Ex_singularII}.
Observe that in OCP \eqref{eq:Fish} with \eqref{eq:Fish_quadI} the Mayer term  $\varphi(x) = -\frac{1}{2} x^2$ is dominant compared to the Lagrange cost $\ell$. Hence, the turnpike leaving arcs approach a part of the  $x-u$ plane with higher running cost $\ell$ but lower terminal cost $\varphi$, see Figure \ref{fig:Ex_quadI}.

  \begin{figure}[htb]
 \includegraphics[scale=0.375]{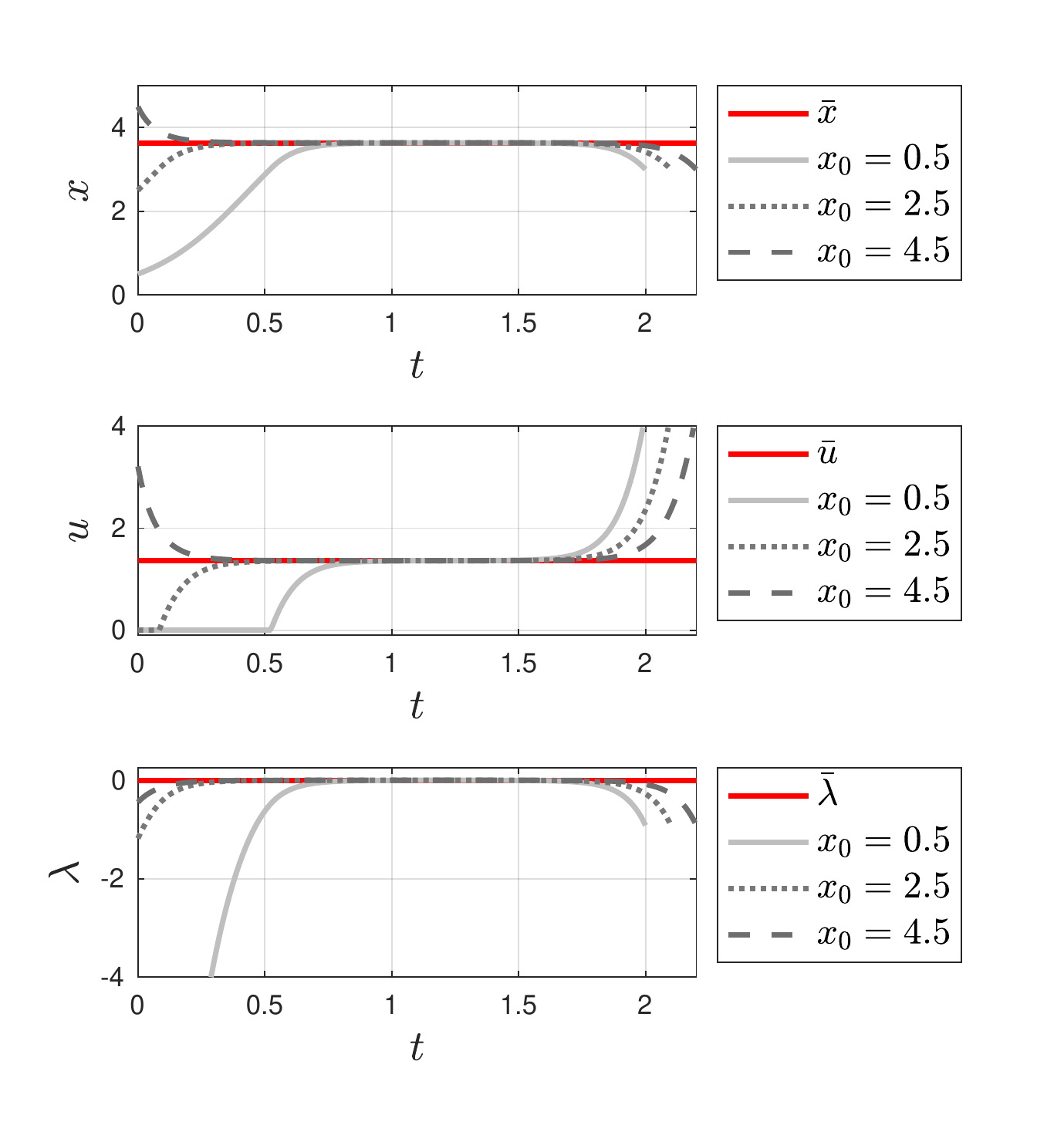}   \includegraphics[scale=0.375]{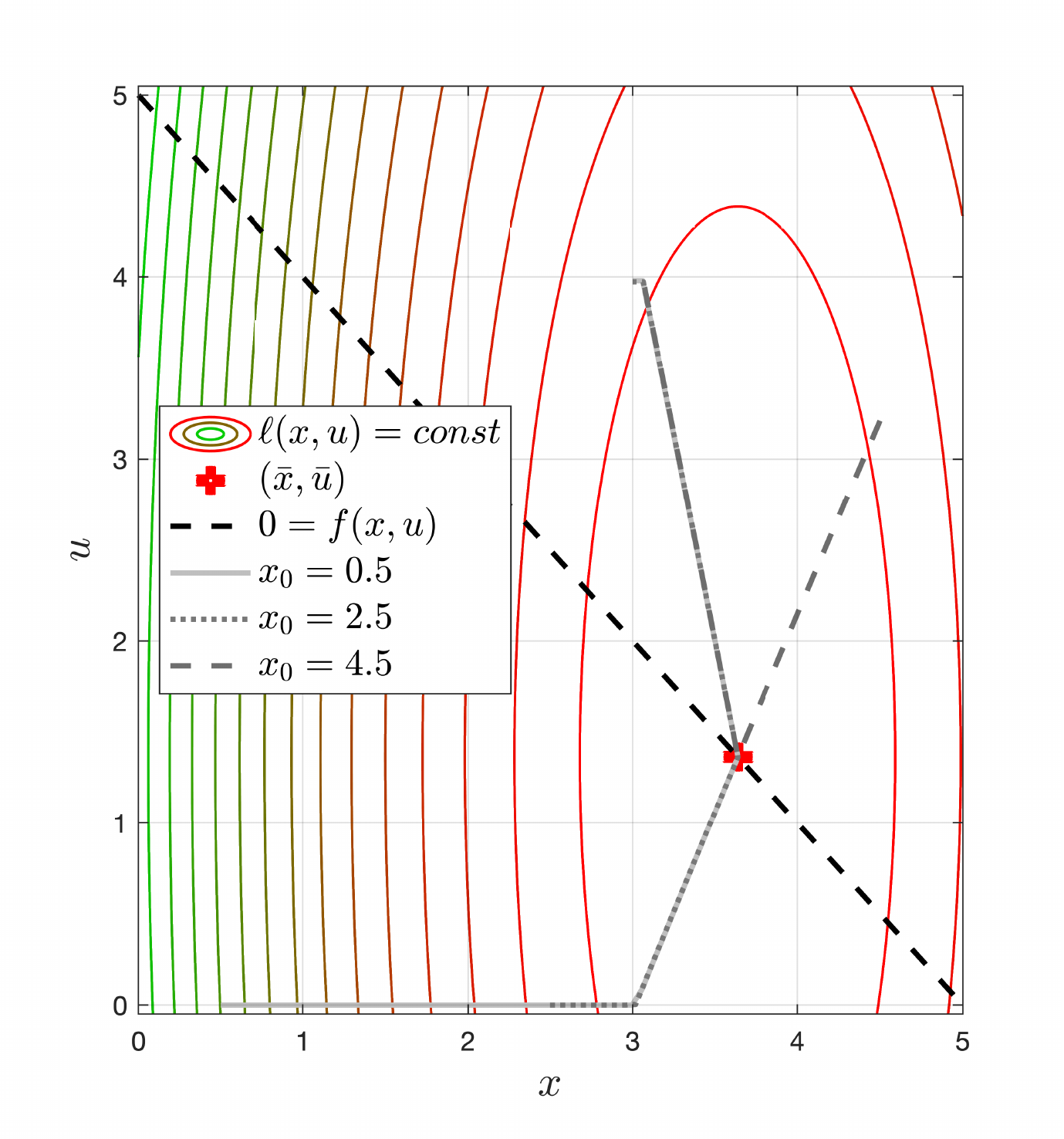}
 \caption{Fish harvest example with quadratic objective and terminal constraint \eqref{eq:Fish_quadII} with pairs of $x_0$ and $T$: $x_0 = 0.5, T = 2.0$;  $x_0 = 2.5, T = 2.1$; and  $x_0 = 4.5, T = 2.2$. Left: optimal state, input and adjoint trajectories. Right: $x-u$ plane with optimal trajectories and contour lines of the Lagrange function $\ell(x,u)$ (coloured), and the set of admissible steady states (dashed black).}\label{fig:Ex_quadII}
 \end{figure} 
 
 Finally, we consider a third variation of OCP \eqref{eq:Fish} by changing \eqref{eq:Fish_quadI} to
 \begin{equation}
 \label{eq:Fish_quadII}
 \ell(x,u) = \dfrac{1}{2}\left( (x - x_c)^2q +  (u - u_c)^2\right), \quad \varphi(x) = 0,\quad \psi(x) = x -x_f, 
 \end{equation}
 with $x_c = 3.64, u_c = 1.37, q = 10$, and $x_f = 3$. The left hand side of Figure \ref{fig:Ex_quadII} shows the numerical solutions for the arising problem with the resulting terminal equality constraint $\psi(x) = x - 3 = 0$.
\begin{table}[htb]
\caption{Comparison of the three variants of the optimal fish harvest problem. } \label{tab:Comparison}
\begin{tabular}{l || c | c}
OCP& Primal Perspective & Adjoint Perspective \\
\hline
\eqref{eq:Fish} & $\exists \tau<T\, \forall t \in [\tau, T]: \ell(x^\star(t), u^\star(t)) < \ell(\bar x, \bar u)$ & $\bar\lambda \not= \lambda^\star(T) =0$ \\
\eqref{eq:Fish} \& \eqref{eq:Fish_quadI}& $\displaystyle \int_0^T \ell(x^\star(t), u^\star(t))\mathrm{d}t + \varphi(x^\star(T))$ & $\bar\lambda \not= \lambda^\star(T) = \varphi_x(x^\star(T))$ \\
\eqref{eq:Fish} \& \eqref{eq:Fish_quadII}& $\displaystyle \int_0^T \ell(x^\star(t), u^\star(t))\mathrm{d}t$ with $\psi(x^\star(T)) = 0$ &
$\bar\lambda \not= \lambda^\star(T)$
\end{tabular}
\end{table}

Observe that, despite the fact that the center point of the convex quadratic objective $x_c = 3.64, u_c = 1.37$ being an unconstrained minimizer of $\ell$ and a steady state, optimal solutions still show the turnpike leaving arc. This is because the problem data \eqref{eq:Fish_quadII} implies that the optimal solutions have to depart from $(\bar x, \bar u)$  to meet the terminal equality constraint $\psi(x(T)) = 0$.

Based on the three considered variants of the optimal fish harvest problem, one can identify three different mechanisms which induce the characteristic turnpike leaving arc, i.e.\ the departure from the neighborhood of the optimal steady state towards the end of the horizon. The middle column of~Table \ref{tab:Comparison} concisely recalls the primal perspective, which we already discussed. Therein, $\tau<T$ denotes a time after which the running cost along the leaving arc is smaller than at the turnpike.

While the middle column describes the turnpike phenomenon that we see in the primal variables, the right hand side column of Table \ref{tab:Comparison} shows the dual perspective, which provides additional valuable insights. Indeed, besides the different Lagrange functions $\ell$, the three considered OCPs differ in terms of the transversality conditions, i.e.\ the terminal constraint for the adjoint variable $\lambda$.

 In case of the original OCP \eqref{eq:Fish}, the absence of a Mayer term implies $\lambda^\star(T) =0$, which differs from the steady value of $\bar\lambda$ given in \eqref{eq:Ex_sing_arc}.
For \eqref{eq:Fish} with \eqref{eq:Fish_quadI}, the Mayer term implies $\lambda^\star(T) = \varphi_x(x^\star(T))$. Here we have chosen the Mayer term such that $\bar\lambda \not= \lambda^\star(T) = \varphi_x(x^\star(T))$. In other words, for the first two variants of the fish harvest problem, it is the difference between the turnpike value of the adjoint and the corresponding boundary conditions---$\bar\lambda \not= \lambda^\star(T)$---which creates the leaving arc. 
The third variant---i.e., \eqref{eq:Fish} with \eqref{eq:Fish_quadII}---includes a terminal equality constraint $\psi(x^\star(T)) = 0$. As we consider $\psi(x) = x - 3$, and $\bar x = 3.64 \not = x^\star(T)$, in this case it is the primal terminal constraint, which forces the optimal solutions to leave the optimal steady state $(\bar x, \bar u)$.

\vspace*{2mm}
\noindent In summary, we arrive at the following observations:
\begin{itemize}
\item The  turnpike entry  requires asymptotic reachability properties. It is triggered by the fact for sufficiently long horizons it is optimal in a time averaged sense to approach the turnpike. 
\item The turnpike as such is the globally optimal steady state of the dynamics under the considered Lagrange function $\ell$.
\item The characteristic leaving arc is induced whenever the terminal constraints (either for the primal or the dual variables) present in the optimality system \eqref{eq:NCO_ct} do not correspond to the turnpike variables. Specifically, this implies that either a primal mismatch between the turnpike state and the terminal equality constraints $\bar x \not \in \{x\in\mbb{R}^{n_x}\,|\, \psi(x) = 0 \}$ or a dual mismatch between $\bar\lambda \not= \lambda^\star(T)$ occurs. In the latter case, it is well-known that the transversality conditions for $\lambda^\star(T)$ \eqref{eq:NCO_dt_trans} / \eqref{eq:NCO_ct_trans} are governed by the interplay of the Mayer term $\varphi$ and the terminal inequality constraints $ \psi(x) \leq 0$, see \citep{Bryson69a}. 
\end{itemize}

\begin{remark}[Exactness of turnpikes]\label{rem:exactTP}
The careful reader will have noticed that in case of  OCP \eqref{eq:Fish} the turnpike property is indeed exact as the optimal solutions enter the optimal steady state, which coincides with the singular arc \eqref{eq:Ex_sing_arc}. Intuitively, it is clear that whenever a turnpike property is exact this induces additional structure to the optimality system \eqref{eq:NCO_ct}. Indeed under fairly mild technical assumptions it has been shown in \cite[Thm. 3]{kit:faulwasser17a} that whenever the only singular arc is a steady state of \eqref{eq:NCO_ct}, then, if they exist, turnpikes are exact. More generally, whenever the set of control values taken is finite, the turnpike will be exact with respect to the inputs. This occurs frequently in case of mixed-integer OCPs, see \cite{tudo:faulwasser20f} for first steps in this direction.
We remark that recently it has been shown $L_1$ tracking terms can also induce exact turnpikes phenomena in finite dimensional OCPs \citep{Gugat20a} and that certain objectives induce exactness of turnpikes also in Hilbert spaces. 

Finally, in \cite[Remark 3 and Lemma 3]{kit:zanon18a} we argued that turnpikes in the interior of $\mbb{Z}$ are never exact provided that the Jacobian linearization at $(\bar x, \bar u)$ is controllable and the steady state problem satisfies differentiability assumptions and LICQ, see also \cite{kit:faulwasser18e_2} for the discrete-time analysis. 
\end{remark}



\subsection{Necessary and Sufficient Conditions for Turnpikes} 
\label{sec:DI_TP} 

\def\N{\mathbb{N}}
\def\R{\mathbb{R}}
\def\U{\mathbb{U}}
\def\X{\mathbb{X}}
\def\Z{\mathbb{Z}}
\def\KK{\mathcal{K}}
\def\NN{\mathcal{N}}
\def\ue{\bar u}
\def\xe{\bar x}
\def\eps{\varepsilon}

As the examples in the previous section showed, the turnpike property is easily observed numerically by looking at simulated optimal trajectories. However, the fact that closed-form expressions of optimal trajectories are hardly ever available indicates that this is not a feasible way to obtain \emph{rigorous} statements about the presence of the turnpike property. 

To this end, sufficient and necessary conditions have been developed in the literature, which we survey in this section. One obvious necessary condition is approximate asymptotic controllability to $\bar \xi$. This property demands that for all neighborhoods $\NN$ of $\bar \xi$ there exists $\tau>0$ such that for all initial conditions $x_0\in\X_0$ and all $T>0$ there is an admissible control $u(\cdot)$ with $\xi(t;x_0)\in \NN$ for all $t\in [\tau,T]$. For affine linear systems $\dot x = Ax+Bu+c$ and $(\xe,\ue)\in{\rm int}\,\Z$, one easily sees that this property is equivalent to stabilizability of the pair $(A,B)$. 

The existing sufficient conditions can be divided into two classes: (i) conditions based on necessary optimality conditions and (ii) conditions based on dissipativity. Regarding (i), \cite{TreZ15} provide a prototype example for finite dimensional systems. Essentially, suitable conditions on the Hamiltonian \eqref{eq:H} of the OCP are formulated, which allow to conclude the continuous-time primal-dual exponential turnpike property. This approach can be extended to larger classes of systems, for instance, to optimal control problems governed by PDEs, as in, e.g., \citep{PorZ13,porretta2016remarks,Zuaz17,GrSS20,Pighin20} or to shape optimization problems \citep{LaTZ20}. It should be mentioned that several of these references consider linear-quadratic problems as a special case or as a local approximation of more general nonlinear problems. This allows to use the Riccati equation as an alternative way to formulate necessary conditions for optimality, and thus also to derive sufficient conditions for turnpike behaviour, typically in form of detectability conditions.

Regarding (i), dissipativity-based approaches rely on the systems-theoretic notion of dissipativity introduced by \cite{Will72a,Will72b}, typically in a strict form. This property requires the existence of so-called storage and supply functions $S:\X\to\R$ and $w:\X\times\U\to\R$, which satisfy the inequality
\begin{align} \nabla S^\top(x) f(x,u) \leq w(x, u)  - \alpha(\|x-\bar x\|+\|u-\bar u\|)\label{eq:diss_ct}\end{align} 
for continuous-time problems and
\begin{align} S(f(x,u)) -S(x) \leq w(x, u) - \alpha(\|x-\bar x\|+\|u-\bar u\|)\label{eq:diss_dt}\end{align} 
for discrete-time problems, where\footnote{As usual, $\KK_\infty$ denotes the space of functions $\alpha:[0,\infty)\to[0,\infty)$ that are continuous, strictly increasing and unbounded with $\alpha(0)=0$.} $\alpha\in\KK_\infty$ and $S$ is bounded from below. For some purposes it is sufficient that only $\|x-\bar x\|$ appears in the argument of $\alpha$ in \eqref{eq:diss_ct} or \eqref{eq:diss_dt}. The link to the OCPs introduced in the last section is then established by specifying $w(x,u)\doteq\ell(x,u)-\ell(\xe,\ue)$.

The earliest result we are aware of that proves that strict dissipativity (plus an appropriate controllability property) implies the turnpike property can be found in \cite[Theorem 4.2]{CaHL91}. The respective Assumption 4.2 in this reference requires the existence of a linear storage functions $S(x) = \bar p^\top x$. The proof, however, is easily extended to general $S$, see \cite[Theorem 5.6]{Grue13} in discrete time and \cite[Theorem 2]{FKJB17} in continuous time. The proof idea is actually very easy; we sketch the idea in discrete time for state turnpike in a problem with terminal cost $\varphi\equiv 0$: assume that the optimal trajectory $x^\star$ stays outside an $\eps$-ball around $\xe$ for an amount of time instances $K(T)$ that grows unboundedly as $T$ tends to infinity. Then summing up \eqref{eq:diss_dt} along the optimal trajectory yields
\begin{align*} V_T(x_0) & = \sum_{t=0}^{T-1} \ell(x^\star(t),u^\star(t)) \ge T\ell(\xe,\ue) - S(x_0) +S(x^\star(T)) + K(T)\alpha(\eps)\\
& \ge C +  T\ell(\xe,\ue) + K(T)\alpha(\eps),\end{align*}
where $C$ is a constant that is independent of $T$ (here the lower bound on \LG{$S$} is important because $x^\star(T)$ may vary with $T$). 
Now, if in addition the inequality 
\begin{align} V_T(x_0)  \le \widetilde C + T\ell(\xe,\ue)\label{eq:cr}\end{align} 
holds for a constant $\widetilde C$ independent of $T$, then we arrive at a contradiction, implying that $K(T)$ must be bounded as $T\to\infty$, from which the turnpike property can be concluded. Inequality \eqref{eq:cr} can be ensured by controllability or stabilizability properties of the dynamics, see, e.g., Theorem 6.4 in \citep{Grue13}. An assumption of this kind is imposed in all dissipativity based proofs of the turnpike property. Notice that the proof as sketched above requires the dissipation inequalities and \eqref{eq:diss_ct} and  \eqref{eq:diss_dt} to hold only along optimal solutions. Moreover, in the continuous-time case, the integral counterpart to \eqref{eq:diss_ct} suffices.\footnote{{The integral counterpart to  \eqref{eq:diss_ct}  reads
\[
S(x(T)) - S(x_0) \leq \int_0^T \ell(x,u) \mathrm{d}t - \int_0^T\alpha(\|x(t)-\bar x\|+\|u(t)-\bar u\|)\mathrm{d}t. 
\] 
This can also be written as
\[
S(x(T)) - S(x_0) \leq V_T(x_0)- \int_0^T\alpha(\|x(t)-\bar x\|+\|u(t)-\bar u\|)\mathrm{d}t, 
\]
which highlights that the dissipation inequality---specifically the storage function $S$---provides a lower bound on $V_T$, see also \citep{Willems71a,tudo:faulwasser20a}.}
}

As this sketch already indicates, dissipativity-based proofs of the turnpike property typically only deliver measure or cardinality turnpike. Yet, under additional conditions on the problem data one can also obtain exponential turnpike \citep{DGSW14}, which is particularly interesting because it provides a direct way to establish exponential turnpike for nonlinear systems without having to resort to linearization arguments. Strict dissipativity inequalities as in \eqref{eq:diss_ct} or \eqref{eq:diss_dt} yield the input-state turnpike property, while if one uses $\alpha(\|x-\bar x\|)$ instead of $\alpha(\|x-\bar x\|+\|u-\bar u\|)$ in these inequalities one obtains the state turnpike property. As there is no direct connection between strict dissipativity and necessary optimality conditions, there is no direct access to the dual variables in this approach and hence---to the best of our knowledge---no direct way to obtain primal-dual turnpike properties from strict dissipativity. 
However, under suitable assumptions on the optimal control problem, the primal-dual turnpike property follows from the input-state turnpike property \citep{FGHS20}. 

As we will discuss in Section \ref{sec:Numerics}, below, strict dissipativity is a very useful property in the context of model predictive control, which allows for deriving qualitative properties of the approximately optimal trajectory generated by this method. It is thus natural to ask how much stronger strict dissipativity is compared to the mere occurrence of the turnpike property. In \citep{GruM16} it was shown for discrete-time controllable systems, that strict dissipativity (with $\alpha(\|x-\bar x\|)$ instead of $\alpha(\|x-\bar x\|+\|u-\bar u\|)$) is equivalent to a robust variant of the state-turnpike property. This robust variant demands that not only optimal but also near-optimal trajectories exhibit turnpike behavior, i.e.\ the turnpike property is robust with respect to perturbations of optimal trajectories. These results can be refined for linear-quadratic problems in both discrete and continuous time, see \citep{GruG18,GruG20}. 

\subsection{Finite and Infinite-Horizon Cases} \label{sec:infHor}

So far we have considered turnpike properties on finite time horizons $T$. While it is important to allow for arbitrarily long horizons $T$, until now each $T$ under consideration was finite. It is, however, also possible (and not uncommon in the literature) to define the turnpike property for optimal control problems on an infinite time horizon. In the discrete-time case, the optimization criterion \eqref{eq:ocp_dt_obj} then becomes 
\begin{equation} V_\infty(x_0) \doteq\inf_{u(\cdot)}  \sum_{k=0}^{\infty} \ell(x(k),u(k)), \label{eq:ocp_dt_inf}\end{equation}
the conditions \eqref{eq:ocp_dt_dyn} and \eqref{eq:ocp_dt_pc} are required for all $k\in \N_0$, and the terminal constraint \eqref{eq:ocp_dt_tc} is dropped. In continuous time, the changes are analogous. 

Definition \eqref{eq:ocp_dt_inf} must be interpreted with care, because in general the infinite sum need not converge. In case an optimal equilibrium exists, the usual remedy is to replace $\ell(x,u)$ by $\ell(x,u)-\ell(\xe,\ue)$. This implies that $V_\infty(\xe)<\infty$ and then additional regularity assumptions are imposed (e.g., bounds and continuity assumptions on the storage function $S$ and on $V_\infty$) in order to ensure that the problem is well posed. This procedure is adopted, e.g., in \cite[Chapter 7]{GruP17}. Alternatively, optimality can be defined using the notion of overtaking optimality \citep{CaHL91,BloH14}, which can in particular be exploited for analyzing turnpike properties for time-varying systems \citep{GrPS18}.

Note that the optimal trajectories are now defined for all $k\in\N_0$. The translation of the turnpike properties from finite to infinite time horizon can then be performed almost formally: it suffices to replace ``$T$'' by ``$\infty$'' in the definitions. For instance, the definition of cardinality turnpike in the discrete-time case becomes:

\begin{definition}[Discrete-time (cardinality) turnpike for infinite horizon problem]
OCP~\eqref{eq:ocp_dt_inf} is said to have a turnpike at $\bar \xi, \xi\in \{x, (x,u), (x,u,\lambda, \mu)\}$ if there exist a function $\nu_\xi:[0,\infty)\to[0,\infty]$ such that, for all $x_0 \in \mbb{X}_0$,
	\begin{equation}\label{eq:TP_dt_inf}
	\#\mathbb{Q}_{\xi,\infty}(\varepsilon)  \leq \nu_\xi(\varepsilon)<\infty\quad \forall\: \varepsilon >0,
	\end{equation}
 where 
	\begin{equation}\label{eq:SetQ_inf}
	\mathbb{Q}_{\xi,\infty}(\varepsilon)\doteq\{t\in\N\hspace*{1mm}|\hspace{1mm} \|\xi^\star(t;x_0)-\bar{\xi}\| > \varepsilon \}
	\end{equation}
	and $\#\mbb{Q}_{\xi,\infty}(\varepsilon)$ is the cardinality of $\mbb{Q}_{\xi,\infty}(\varepsilon)$.
	\begin{itemize}
	\item If $\xi = x$, then OCP~\eqref{eq:ocp_dt_inf} is said to have a \textit{state turnpike} at $\bar x$.
	\item If $\xi = (x, u)$, then OCP~\eqref{eq:ocp_dt_inf} is said to have an \textit{input-state turnpike} at $(\bar x, \bar u)$.
	\item If $\xi = (x, u, \lambda, \mu)$, then OCP~\eqref{eq:ocp_dt_inf} is said to have a \textit{primal-dual turnpike} at $(\bar x, \bar u, \bar 
	\lambda, \bar \mu)$.
	\end{itemize}
\end{definition}
It is easy to check that this definition implies that for each $x_0\in\X_0$ and each $\eps>0$ there is $T_\eps>0$ such that 
\[ \|\xi^\star(t;x_0)-\bar{\xi}\| \le  \varepsilon \]
for all $t\ge T_\eps$---which is nothing else than saying that $\xi^\star(t;x_0)$ converges to $\bar{\xi}$ as $t\to\infty$. The infinite-horizon cardinality turnpike property hence demands that all optimal trajectories (possibly including the corresponding inputs and costates) converges towards the same steady state.

Interestingly, the analogous infinite-horizon variant of the continuous-time measure turnpike property  (which is obtained replacing $T$ by $\infty$ in \eqref{eq:SetTheta} and \eqref{eq:TP_ct}) does not imply this convergence. This is because the mere fact that the measure $\nu_\xi(\eps)$ is finite does not exclude the existence of arbitrary large times $t\in\Theta_{\xi,\infty}$. However, for instance for the state trajectory, convergence can be recovered if the the state trajectory is uniformly continuous in $t$, because the Lebesgue measure of $\Theta_{\xi,\infty}\cap[t_1,\infty)$ must tend to $0$ as $t_1$ tends to $\infty$. This means that the length of time intervals contained in $\Theta_{\xi,\infty}\cap[t_1,\infty)$ tends to zero, and uniform continuity then implies that the deviation of the trajectories from the turnpike during these time intervals also becomes arbitrarily small as $t_1$ tends to $\infty$. 
Alternatively, \citep[Chap.\ 4.6]{CaHL91} discuss an approach to guarantee primal-dual stability of infinite-horizon continuous-time optimal solutions via the Lyapunov function candidate $\lambda^\top x$, certain convexity assumptions and LaSalle's invariance principle. Moreover, \cite[Thm.\ 2 and Thm.\ 4]{tudo:faulwasser20a} show convergence and local exponential stability of primal solutions via a dissipativity assumption.

In the case of the exponential turnpike, the infinite-horizon version is obtained replacing \eqref{eq:expTP_dt} by 
\[ \|\xi^\star(t;x_0)-\bar{\xi}\| \leq C\rho^{t}, \]
which is required for all $t\ge 0$. This implies exponential convergence of the optimal trajectories towards the turnpike $\bar\xi$ in both discrete and continuous time. In case the stronger estimate $\|\xi^\star(t;x_0)-\bar{\xi}\| \leq C\|x_0-\xe\|\rho^{t}$ holds, we even obtain exponential stability of the optimally controlled system.  Conversely, for infinite-horizon OCPs, one can show as a local result that exponential stability of the infinite-horizon dual variables implies strict input-state dissipativity \citep{tudo:faulwasser20a}.

The generating mechanisms for infinite-horizon turnpike behaviour are very similar to those in the finite-horizon case. Under suitable regularity conditions on the optimal control problem, it can even be rigorously shown that finite-horizon turnpike behaviour occurs if and only if infinite-horizon turnpike behaviour occcurs. Conditions under which this is true are given, e.g., in \cite[Theorem 3.1.4]{Zasl06} or in \cite[Theorem 3]{Zasl08} in continuous time, and in \citep{GrKW17} in discrete time (see also \cite[Theorem 2.2]{Zasl14} for conditions under which the infinite-horizon turnpike property implies the finite-horizon property in discrete time). To give an impression about the flavor of such conditions, we briefly explain the two conditions needed in \cite[Theorem 3.1]{GrKW17}: The first condition is that the optimal value functions for finite and infinite horizon are bounded on bounded sets by bounds that are independent of the horizon length. The second condition demands that trajectories with bounded objective remain in bounded sets, again uniformly, i.e.\ the size of the bounded set is determined by the bound on the objective and does not depend on the time horizon. Here it should be mentioned that the robust definition of the turnpike property discussed in the last pararaph of Section \ref{sec:DI_TP} is needed for this result.

\begin{remark}[Relation between $\bar\lambda$ and $V_\infty$]
In case of a primal-dual turnpike as per Definition \ref{def:MeasTP} one immediately has that for sufficiently large horizons $T$, there exists $\tau \in [0,T]\setminus\Theta_{(x,u,\lambda,\mu), T}(\varepsilon)$ such that
\[
\textstyle\nabla V_T(x(\tau)) =\lambda^\star(\tau) \approx \bar\lambda,
\]
i.e. close to the turnpike $\bar\lambda$ is an approximation of $\nabla V_T$.
Under appropriate regularity assumptions we can also make statements about the dual variables in the infinite-horizon case. For instance, under suitable differentiability assumptions it can be shown that the gradient of the optimal value function of the infinite-horizon OCP, $V_\infty$,  satisfies
\[
\nabla V_\infty(\bar x) = \bar\lambda = -\nabla S(\bar x),
\]
cf. \cite[Thm. 10]{tudo:faulwasser20a}.

 In view of Table \ref{tab:Comparison} it also appears natural to conjecture that for $T=\infty$ the optimal adjoint trajectory $\lambda^\star$ approaches its turnpike value. Indeed, a recent analysis in \citep{tudo:faulwasser20a} shows that assuming dissipativity one can provide a positive answer to this classical problem of adjoint boundary conditions for infinite-horizon OCPs posed by \cite{Halkin74a}.
\end{remark}


\section{Exploitation of Turnpikes in Numerics and Receding-horizon Control} \label{sec:Numerics}

There are at least three different ways how the turnpike property can be exploited for the numerical computation of optimal trajectories: (i) splitting the optimization horizon at the turnpike, (ii) receding-horizon approximation, and (iii) exploiting the turnpike phenomenon in the numerical discretization appearing in (ii).
  All of these approaches are particularly efficient when optimal trajectories on long, possible infinitely long time intervals $[0,T]$, $T\in\R\cup\{\infty\}$, are sought. 
We note that when the length of the time horizon $T$ increases, the problem size---e.g., the number of decision variables in direct numerical solution methods---grows while the conditioning of the problem often deteriorates. At the same time, the effort to solve underlying linear systems of equations grows super-linearly with their
dimension. Put differently, we have the following empirical insight for OCPs more complex than a linear-quadratic regulator:  the longer the optimization horizon the more difficult it is
to solve an OCP numerically. Moreover, the direct computation of infinite-horizon optimal trajectories is challenging and often only possible using dynamic programming techniques, which become numerically infeasible already for moderate state dimension, see, e.g., \citep{Bert95}. This motivates techniques that circumvent the direct computation of optimal trajectories on long and infinite horizons.

The \emph{first approach} is the most straightforward one: Given the turnpike equilibrium $\xe$, the optimization horizon $T>0$ (possibly $T=\infty$), an initial condition $x_0$ and, in case $T<\infty$, a terminal condition $x(T)\in\X_f\doteq \{x\in\mbb{X}\,|\, \psi(x) \leq 0\}$, we compute 
\begin{itemize}
\item an optimal trajectory $x_1(\cdot)$ with finite horizon $T_1<T$ and initial and terminal conditions $x_1(0)=x_0$, $x_1(T_1)=\xe$
\item in case $T<\infty$, an optimal trajectory $x_2(\cdot)$ with horizon $T_2<T-T_1$ and initial and terminal conditions $x_2(0)=\xe$, $x_2(T_2)\in\X_f$.
\end{itemize}
An approximation of the optimal trajectory is then obtained by
\[ \tilde x(t) = \left\{ \begin{array}{ll} x_1(t), & t\in[0,T_1]\\
x^e, & t\in [T_1,T-T_2]\\
x_2(t-T+T_2), & t \in [T-T_2,T]
\end{array}\right.
\]
in case $T<\infty$ and
\[ \tilde x(t) = \left\{ \begin{array}{ll} x_1(t), & t\in[0,T_1]\\
x^e, & t\in [T_1,\infty)
\end{array}\right.
\]
in case $T=\infty$. The idea for this construction goes back at least as far as \citep{Wilde72a}, see also \citep{AndK87}, in which also the error of the method is analyzed for particular classes of optimal control problems. This construction has also been suggested for PDEs, see~\cite[Rem. 7]{hernandez2019greedy}. If the turnpike is not exact, the method causes an error because the approaching arc---which is approximated by $x_1(\cdot)$---does not reach the turnpike equilibrium $\xe$ exactly but only a neighbourhood of $\xe$. Likewise, the middle part of the optimal solution is only near $\xe$ but does not exactly coincide with $\xe$. Finally, the leaving arc, which is approximated by $x_2(\cdot)$, does in general not start exactly at $\xe$ but only near $\xe$. The resulting error can be estimated if the speed of convergence towards the turnpike is known (as in the case of exponential turnpike) and if the sensitivity of the optimal value with respect to the initial condition $x_0\approx \xe$ is independent of $T$, i.e.\ if the optimal value for initial conditions $x_0$ satisfying $\|x_0-\xe\|\le\delta$ satisfies $|V_T(x_0)-V_T(\xe)|\le \eps$ with $\eps>0$ independent of $T$ and $\eps\to 0$ as $\delta\to 0$. In this case the error of each of the trajectory pieces constituting $\tilde x(\cdot)$ and thus the overall error can be estimated. 

The \emph{second approach} uses Model Predictive Control (MPC), also known as Receding Horizon Control. In this approach, we fix a finite horizon $T_{opt}<T$ and a control horizon $\delta<T_{opt}$, set $t:=0$, $x_{MPC}(0):=x_0$, and proceed as follows:
\begin{enumerate}
\item \label{step1}Solve the optimal control problem with horizon $T_{opt}$ (possibly with additional or modified terminal conditions) with initial condition $x_0 = x_{MPC}(t)$ and denote the numerical approximation of the optimal solution by $x_{opt}(\cdot)$
\item set $x_{MPC}(t+\tau) := x_{opt}(\tau)$ for $\tau\in[0,\delta]$, set $t:=t+\delta$
\item if $T<\infty$, set $T_{opt}:=\min\{T_{opt},T-t\}$, $\delta:=\min\{\delta,T-t\}$
\item if $T_{opt}>0$, go to {\sf \textbf{1.}}
\end{enumerate}

Compared to the first approach, MPC has the advantage that the optimization is carried out repeatedly. Efficient modern algorithms can perform one of the optimizations in Step {\sf \textbf{1.}}\ of the MPC algorithm in the range of milliseconds or even below \citep{HoFD11}. This means that for many applications the approach is real-time capable and thus provides the control in form of a feedback law that is able to react to perturbations or modelling errors. Yet, even if for a specific OCP the MPC methodology is not real-time capable, e.g., due to the dimensionality of the problem or due to the complexity of the dynamics, MPC is a powerful approach to reduce the computational complexity of the problem and allows to approximate problems that cannot be solved directly. As such, MPC can be seen as a complexity reduction technique in time.

Under essentially the same conditions as for the first approach it can be shown that this method yields approximately optimal trajectories \citep{Grue16}. If one additionally assumes strict dissipativity \eqref{eq:diss_ct} or \eqref{eq:diss_dt}, it can also be shown that $x_{MPC}$ has the turnpike property, at least aproximately \citep{GruS14}. Suitable terminal conditions on the optimal control problem to be solved in Step {\sf \textbf{1.}}\ can improve both the performance properties as well as the qualitative behavior of $x_{MPC}$ \cite[Chapter 7]{kit:zanon18a,FaGM18,GruP17}. Interestingly, so far this analysis appears to be limited to the case $T=\infty$ as we are not aware of a rigorous analysis for the case of a finite $T$. Yet, we expect that similar arguments work for finite $T$. 

The \emph{third approach} is intertwined with the second approach, as it concerns efficient numerical discretization for computing the solutions of the optimal control problem in Step {\sf \textbf{1.}}\ of the MPC algorithm, particularly for problems governed by partial differential equations (PDEs). It relies on the obvious fact that only $x_{opt}|_{[0,\delta]}$ is needed for generating $x_{MPC}$ in Step {\sf \textbf{2.}} of the MPC algorithm, hence there is no need to compute $x_{opt}(\tau)$ accurately for $\tau\in(\delta, T_{opt}]$. 

To this end, the following interpretation of the turnpike property is useful: it implies that optimal trajectories are essentially independent of the initial condition and---if present---of the terminal condition, for times that are not close to the initial or terminal time, respectively. In the particular case of the exponential turnpike property, one can even quantify that the influence of initial and terminal condition is damped exponentially as $t$ increases or decreases, respectively. For several classes of PDE governed optimal control problems it turns out \citep{GrSS19,GrSS20} that the same mechanism that generates the exponential turnpike property also implies that perturbations acting on the system at a certain point $t_{pert}$ in time are damped out exponentially, i.e.\ their influence on the solution decreases exponentially in $|t-t_{pert}|$. This property can be called an \emph{exponential sensitivity} property and was also established recently for other types of optimal control problems, see, e.g., \cite{SFZZ19}.
Since discretization errors can be interpreted as perturbations, this implies that increasingly coarse discretizations towards the end of the horizon do only marginally affect the accuracy of $x_{opt}|_{[0,\delta]}$. This observation forms the basis for the design of efficient adaptive discretizations in time and space, which can be generated using appropriate goal oriented error estimators \citep{GrSS20b}. Moreover, the third approach has also been considered for solving  OCPs concerning chemical processes wherein in the middle of the horizon a less dense discretization scheme was used while no error bounds are given, see \citep{Sahlodin15}.

\begin{figure}[t]
	\centering
	\scalebox{.7}{\input{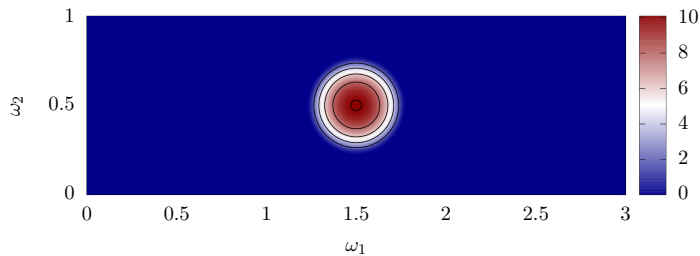}}
	\caption[Constant reference $x_c(\omega)$]{Constant reference trajectory $x_c(\omega)$}
	\label{fig:reference}
\end{figure}

\begin{figure}
	\begin{tikzpicture}[scale=0.95]
	\def\d{2.3}
	\def\t{0.5cm}
	\node (0,0){};
	\node[label={[label distance=0.2cm,text depth=-1ex]above:time $t$}] at (0.5\linewidth,1.4) {};
	\node[label={[label distance=0.2cm,text depth=-1ex]above:goal oriented error estimator}] at (0.2\linewidth,1.4) {};
	\node[label={[label distance=0.2cm,text depth=-1ex]above:standard error estimator}] at (0.8\linewidth,1.4) {};
	\draw [very thick](0.5\linewidth,1.4) -> (0.5\linewidth,-5.5*\d) node [above right] {};
	\draw [arrow,very thick](0.5\linewidth,-6.5*\d) -> (0.5\linewidth,-7.5*\d) node [above right] {};
	\draw [very thick,dotted] (0.2\linewidth,-5.8*\d)-> (0.2\linewidth,-6.2*\d)[]{};
	\draw [very thick,dotted] (0.5\linewidth,-5.8*\d)-> (0.5\linewidth,-6.2*\d)[]{};
	\draw [very thick,dotted] (0.8\linewidth,-5.8*\d)-> (0.8\linewidth,-6.2*\d)[]{};
	
	\draw [very thick](0.5*\linewidth-\t,0) -> (0.5*\linewidth+\t,0);
	\draw [very thick](0.5*\linewidth-\t,0) -> (0.5*\linewidth+\t,0) node [ align=right,above right] {};
	\node[label={[anchor=south east]west:$t=0$}] at (0.5\linewidth,0.2cm) {};
	\draw [very thick](0.5*\linewidth-\t,-1*\d) -> (0.5*\linewidth+\t,-1*\d);;
	\node[label={[anchor=south east]west:$\tau=0.5$}] at (0.5\linewidth,-4.4) {};
	
	\draw [very thick](0.5*\linewidth-\t,-2*\d) -> (0.5*\linewidth+\t,-2*\d) node [ align=right,above right] {};
	\draw [very thick](0.5*\linewidth-\t,-3*\d) -> (0.5*\linewidth+\t,-3*\d) node [ align=right,above right] {};
	\draw [very thick](0.5*\linewidth-\t,-4*\d) -> (0.5*\linewidth+\t,-4*\d) node [ align=right,above right] {};
	\draw [very thick](0.5*\linewidth-\t,-5*\d) -> (0.5*\linewidth+\t,-5*\d) node [ align=right,above right] {};
	\draw [very thick](0.5*\linewidth-\t,-7*\d) -> (0.5*\linewidth+\t,-7*\d);
	\node[label={[anchor=south east]west:$T=10$}] at (0.5\linewidth,-15.9) {};

	\node[inner sep=0pt,anchor=west] (whitehead) at (0,0)
	{\includegraphics[width=0.4\linewidth]{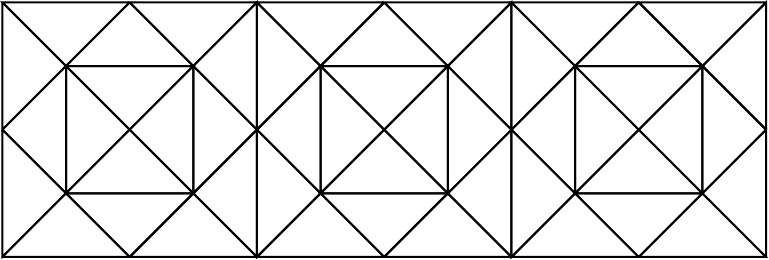}};
	\node[inner sep=0pt,anchor=east] (whitehead) at (0.2\linewidth,0)
	{\includegraphics[width=0.2\linewidth]{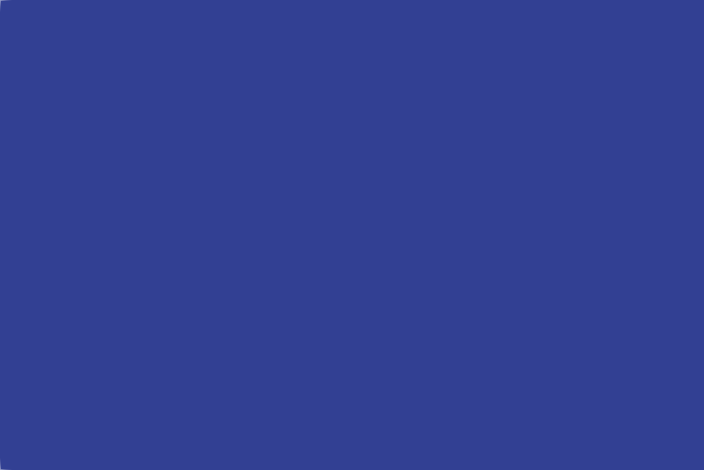}};
	
	\node[inner sep=0pt,anchor=west] (whitehead) at (0,-1*\d)
	{\includegraphics[width=0.4\linewidth]{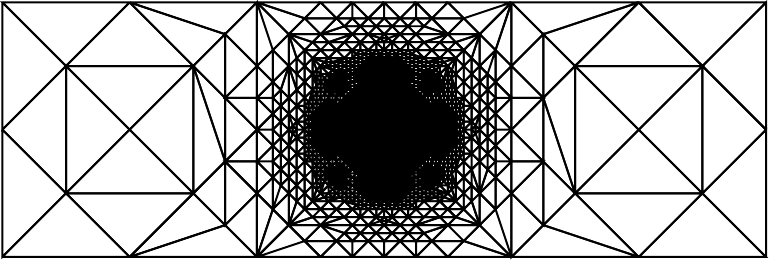}};
	\node[inner sep=0pt,anchor=east] (whitehead) at (0.2\linewidth,-1*\d)
	{\includegraphics[width=0.2\linewidth]{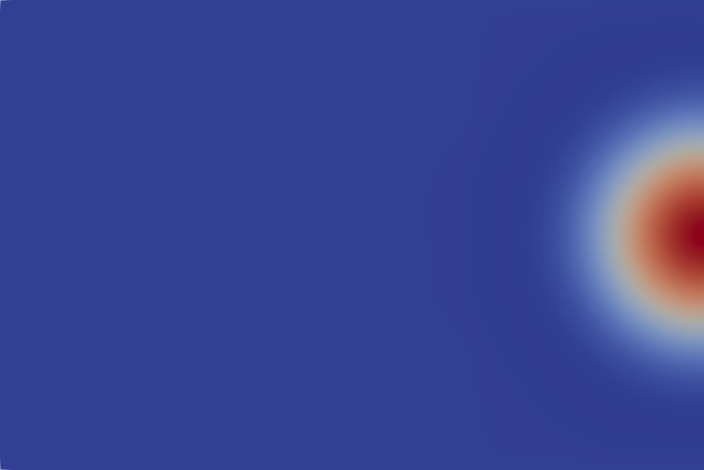}};
	\node[inner sep=0pt,anchor=west] (whitehead) at (0,-2*\d)
	{\includegraphics[width=0.4\linewidth]{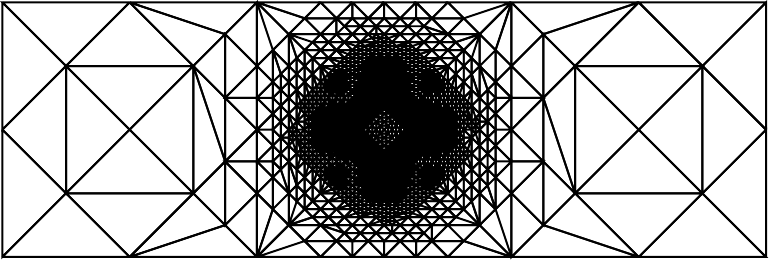}};
	\node[inner sep=0pt,anchor=east] (whitehead) at (0.2\linewidth,-2*\d)
	{\includegraphics[width=0.2\linewidth]{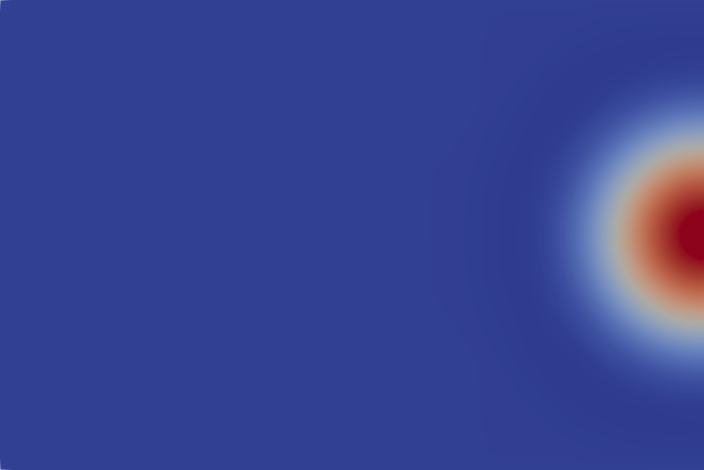}};
	\node[inner sep=0pt,anchor=west] (whitehead) at (0,-3*\d)
	{\includegraphics[width=0.4\linewidth]{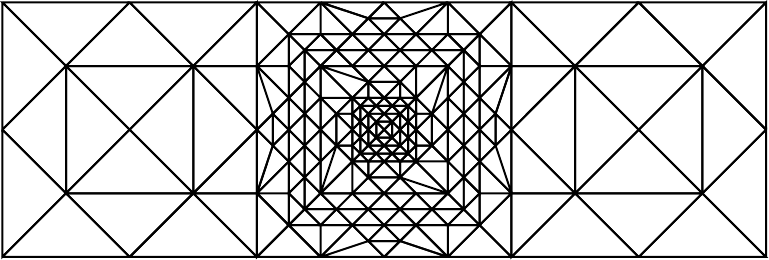}};
	\node[inner sep=0pt,anchor=east] (whitehead) at (0.2\linewidth,-3*\d)
	{\includegraphics[width=0.2\linewidth]{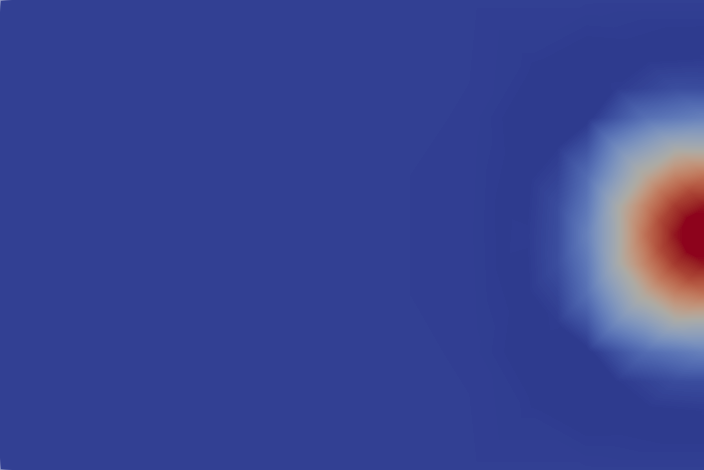}};
	\node[inner sep=0pt,anchor=west] (whitehead) at (0,-4*\d)
	{\includegraphics[width=0.4\linewidth]{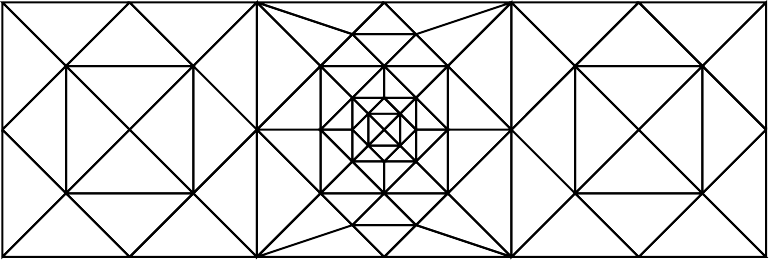}};
	\node[inner sep=0pt,anchor=east] (whitehead) at (0.2\linewidth,-4*\d)
	{\includegraphics[width=0.2\linewidth]{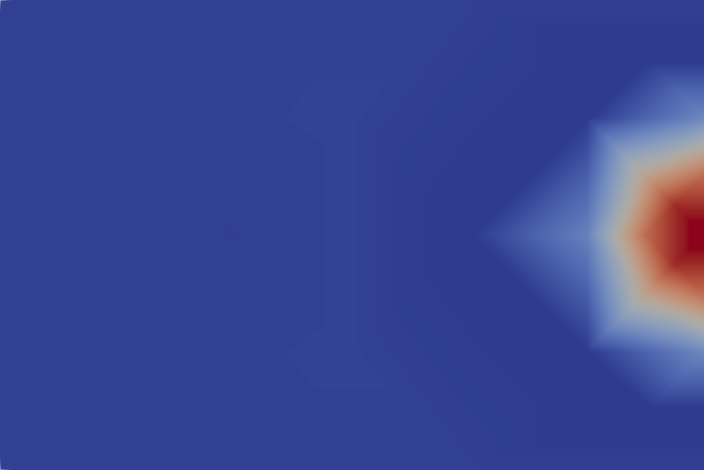}};
	\node[inner sep=0pt,anchor=west] (whitehead) at (0,-5*\d)
	{\includegraphics[width=0.4\linewidth]{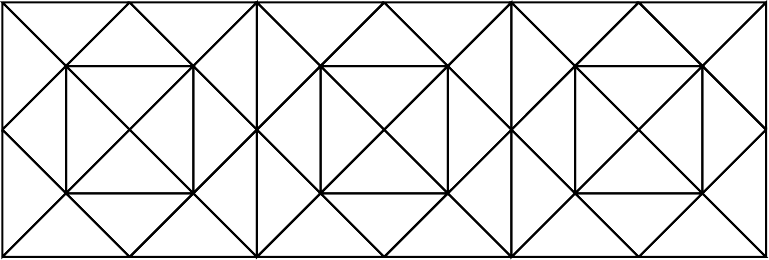}};
	\node[inner sep=0pt,anchor=east] (whitehead) at (0.2\linewidth,-5*\d)
	{\includegraphics[width=0.2\linewidth]{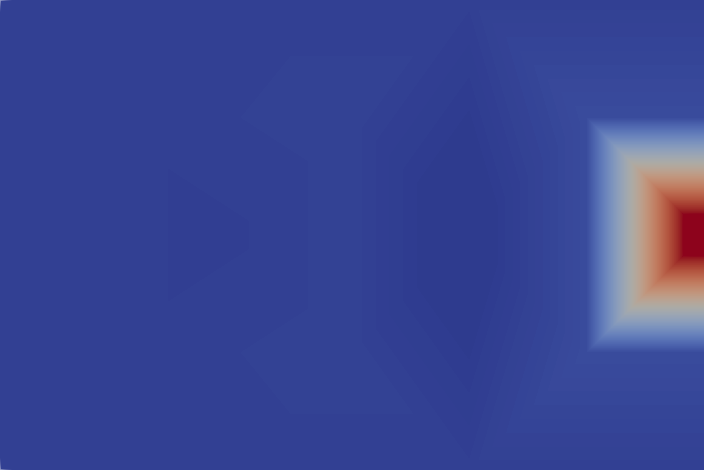}};
	\node[inner sep=0pt,anchor=west] (whitehead) at (0,-7*\d)
	{\includegraphics[width=0.4\linewidth]{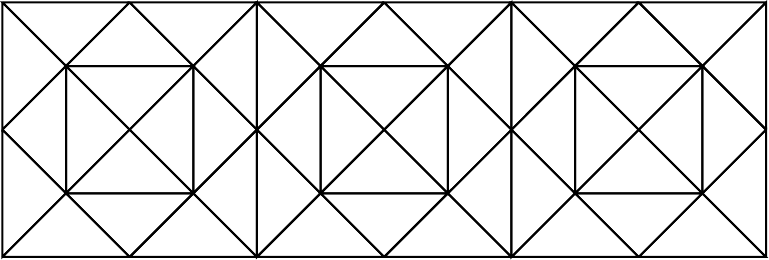}};
	\node[inner sep=0pt,anchor=east] (whitehead) at (0.2\linewidth,-7*\d)
	{\includegraphics[width=0.2\linewidth]{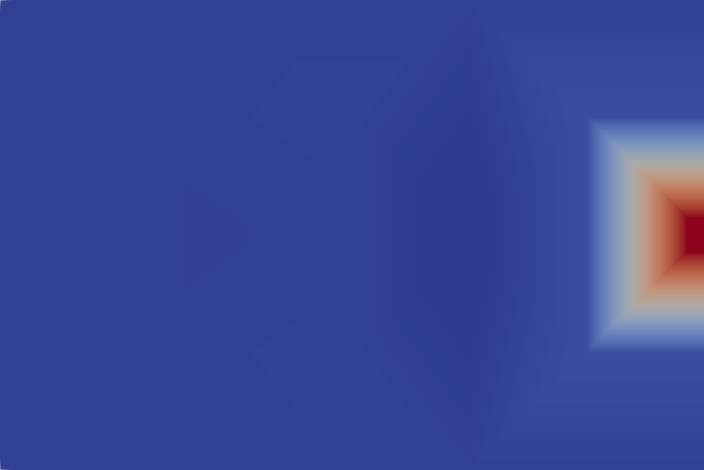}};

	\node[inner sep=0pt,anchor=east] (whitehead) at (\linewidth,0)
	{\includegraphics[width=0.4\linewidth]{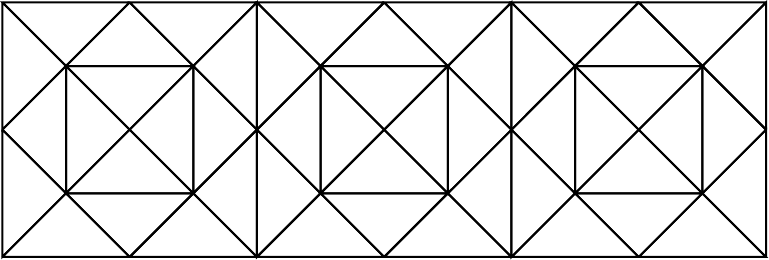}};
	\node[inner sep=0pt,anchor=west] (whitehead) at (0.8\linewidth,0)
	{\includegraphics[width=0.2\linewidth]{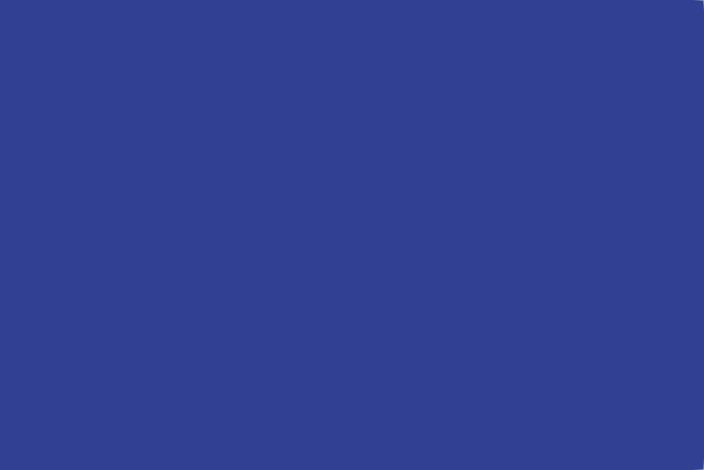}};
	\node[inner sep=0pt,anchor=east] (whitehead) at (\linewidth,-1*\d)
	{\includegraphics[width=0.4\linewidth]{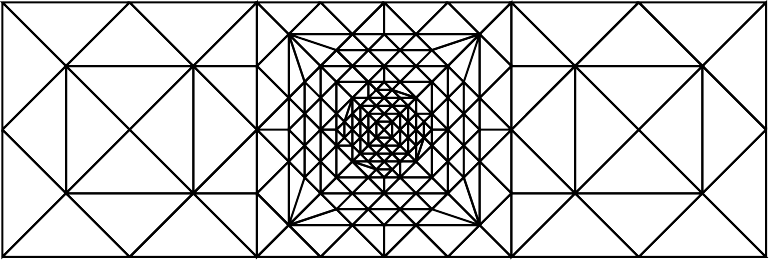}};
	\node[inner sep=0pt,anchor=west] (whitehead) at (0.8\linewidth,-1*\d)
	{\includegraphics[width=0.2\linewidth]{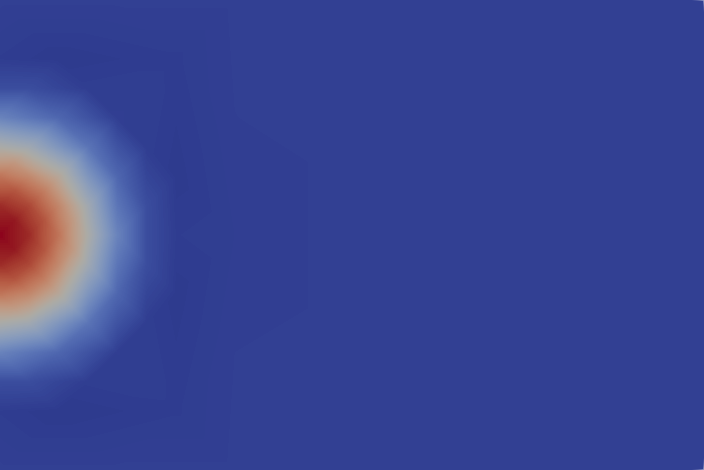}};
	
	\node[inner sep=0pt,anchor=east] (whitehead) at (\linewidth,-2*\d)
	{\includegraphics[width=0.4\linewidth]{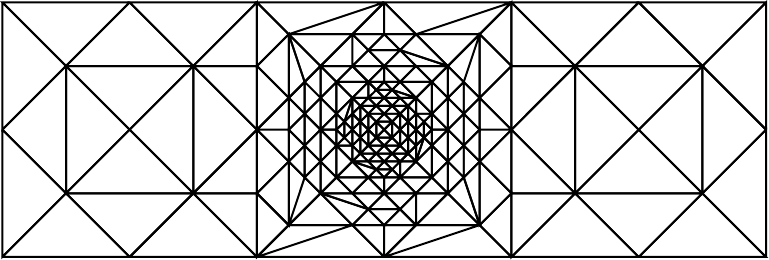}};
	\node[inner sep=0pt,anchor=west] (whitehead) at (0.8\linewidth,-2*\d)
	{\includegraphics[width=0.2\linewidth]{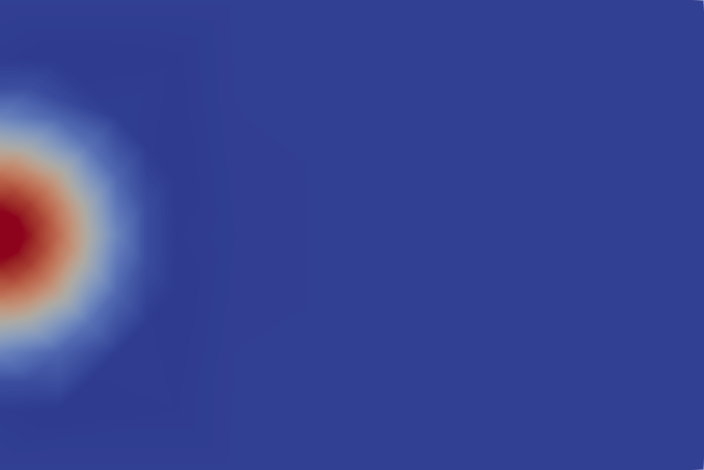}};
	\node[inner sep=0pt,anchor=east] (whitehead) at (\linewidth,-3*\d)
	{\includegraphics[width=0.4\linewidth]{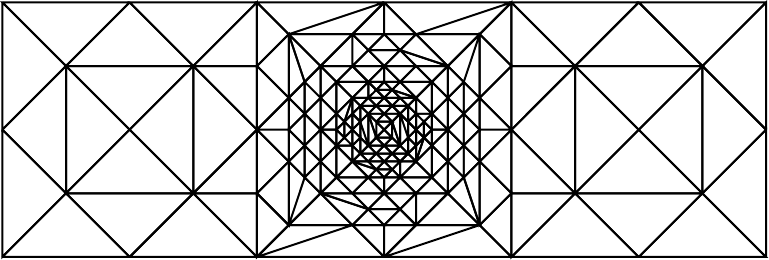}};
	\node[inner sep=0pt,anchor=west] (whitehead) at (0.8\linewidth,-3*\d)
	{\includegraphics[width=0.2\linewidth]{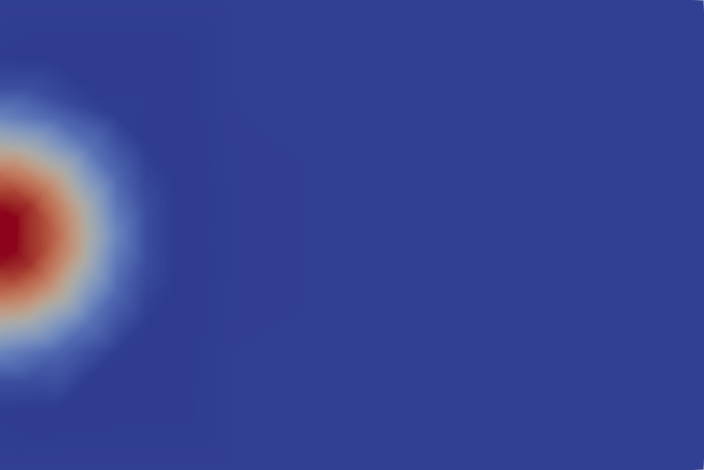}};
	\node[inner sep=0pt,anchor=east] (whitehead) at (\linewidth,-4*\d)
	{\includegraphics[width=0.4\linewidth]{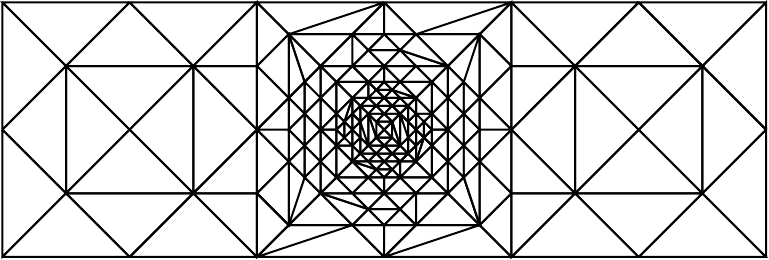}};
	\node[inner sep=0pt,anchor=west] (whitehead) at (0.8\linewidth,-4*\d)
	{\includegraphics[width=0.2\linewidth]{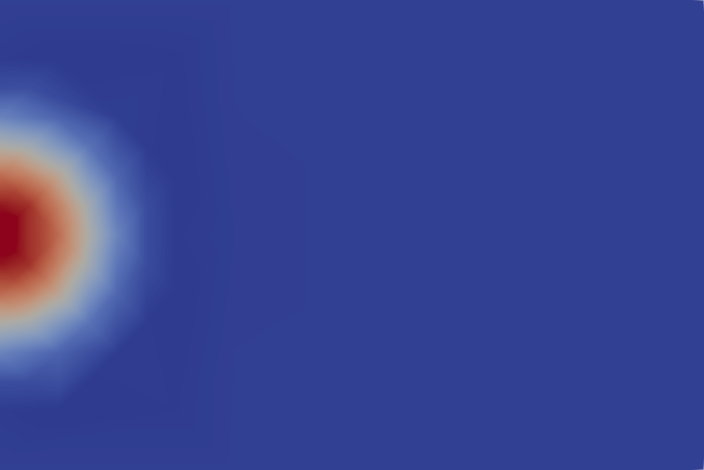}};
	\node[inner sep=0pt,anchor=east] (whitehead) at (\linewidth,-5*\d)
	{\includegraphics[width=0.4\linewidth]{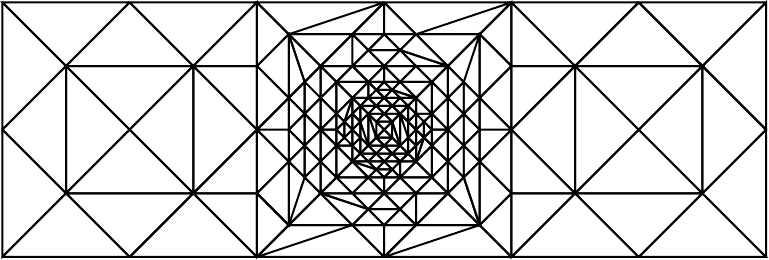}};
	\node[inner sep=0pt,anchor=west] (whitehead) at (0.8\linewidth,-5*\d)
	{\includegraphics[width=0.2\linewidth]{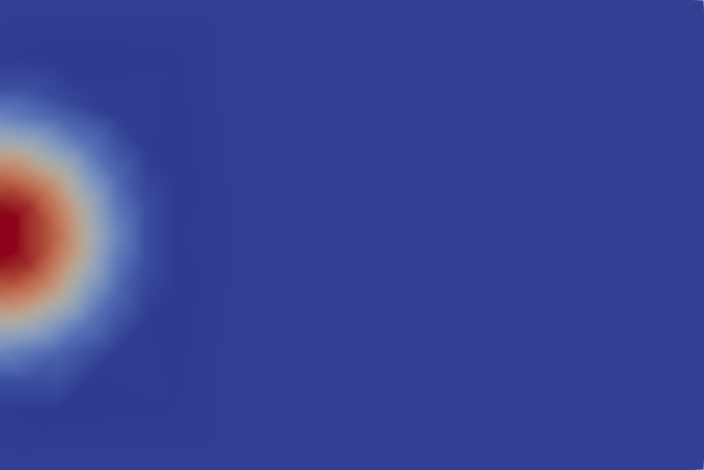}};
	\node[inner sep=0pt,anchor=east] (whitehead) at (\linewidth,-7*\d)
	{\includegraphics[width=0.4\linewidth]{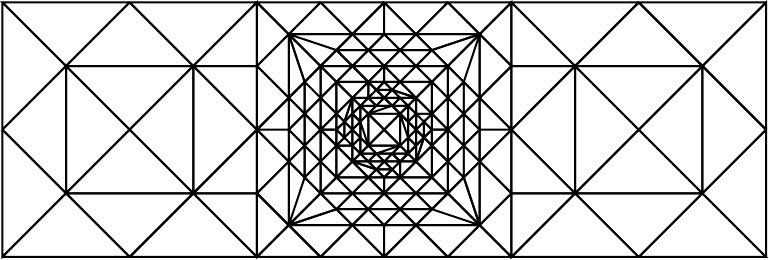}};
	\node[inner sep=0pt,anchor=west] (whitehead) at (0.8\linewidth,-7*\d)
	{\includegraphics[width=0.2\linewidth]{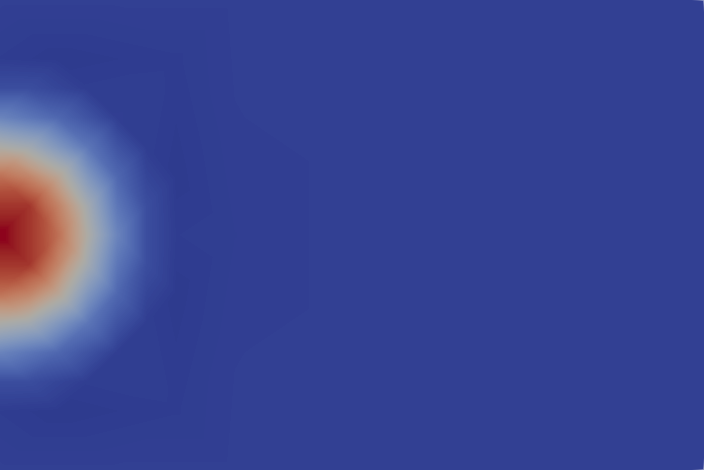}};
	\end{tikzpicture}
	\caption[Evolution of adaptively refined space grids for a linear quadratic problem]{Evolution of adaptively refined space grids for an MPC adapted goal oriented error estimator (left) and a standard error estimator (right) with 5825 and 5478 total space grid points, respectively.}
	\label{fig:grid}
\end{figure}
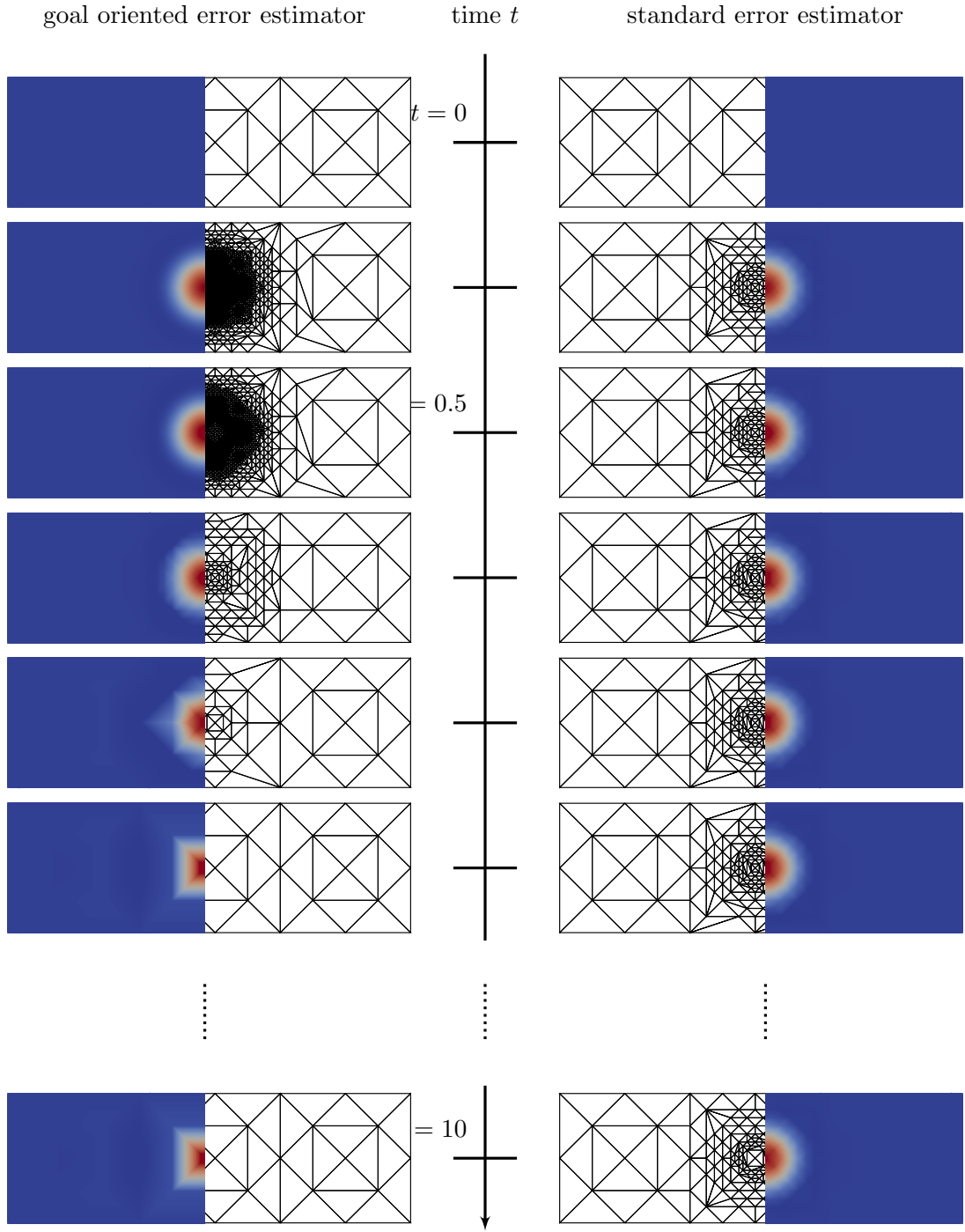

We illustrate the third approach and its use in MPC by means of a numerical example. We consider the linear heat equation with distributed control
\begin{align}
\nonumber
&&\dot{x}(t,\omega) &=   0.1\Delta x(t,\omega) + u  &&\text{for }(t,\omega)\in(0,T) \times \Omega \\
\label{eq:numerics:linear_distributed}
&&x(t,\omega) &= 0 &&\text{for } (t,\omega)\in (0,T)\times \partial\Omega\\
\nonumber
&&x(0,\omega) &= 0 &&\text{for } \omega \in \Omega,
\end{align}
with $\Omega = [0,3] \times [0,1]$ and $T=10$. We consider the cost function 
\[ \ell(x(t,\cdot),u(t,\cdot)) = \frac 12 \|x(t,\cdot)-x_{ref}(t,\cdot)\|^2_{L_2(\Omega)} + \alpha\|u\|_U^2, \]
with $\alpha=10^{-3}$, where $x_{ref}$ may either be the constant reference 
\[ x_{ref}(t,\omega) = x_c(\omega) := g\left(\left\|\omega - \begin{pmatrix}
	1.5\\0.5\end{pmatrix}\right\|\right) \]
depicted in Figure \ref{fig:reference}, or the time-varying reference 
\[ x_{ref}(t,\omega) = x_{tv}d(t,\omega) := g\left(\left\|\omega - 
	\omega(t)\right\|\right), \]
with 
\[ \omega(t) = 	\left( \, 1.5 - \cos\left( \pi\left( \frac{t}{10} \right)\right),\;
\left\vert \cos\left(\pi\left(\frac{t}{10}\right)\right)\right\vert \, \right) ^\top.\]
In both cases we use $g(r) = 10 \exp(1-1/(1-100r^2/9))$ for $s<1$ and $g(r)=0$ else.
It was shown in \citep{GrSS19} that this problem has the exponential turnpike and the exponential sensitivity property. 
Figure \ref{fig:grid} shows the adaptive grid over time for the numerical solution of one optimal control problem in Step~{\sf \textbf{1.}}\ of the MPC algorithm for the static reference. 
On the right, the adaptivity criterion is obtained from a usual error estimator on $[0,T]$, resulting in a grid that is refined in almost all time steps. On the left, a goal oriented error estimator is used, that only aims at providing an accurate solution on the interval $[0,\delta]$ that is actually implemented in the MPC algorithm. here with $\delta=0.5$. The fineness of the grid is chosen such that the total number of space grid points is approximately equal on both sides. Just as the turnpike and sensitivity theory predicts, the accuracy of the solution for $t\in[0,\delta]$ is only marginally affected by errors for large $t$, hence the error estimator allows for coarse grids. 

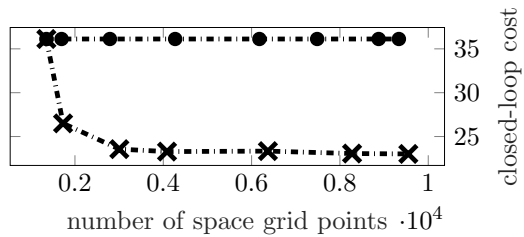
\begin{figure}[htb]
	\centering
	\scalebox{.9}{
%
%
\definecolor{mycolor1}{rgb}{0.00000,0.44700,0.74100}%
\begin{tikzpicture}

%
%
%
%

\begin{axis}[%
width=2.5in,
height=0.8in,
at={(2.7in,1.5in)},
scale only axis,
xlabel style={font=\color{white!15!black}},
xlabel={number of space grid points},
yticklabel pos = right,
ylabel near ticks,
yminorticks=false,
ylabel style={font=\color{white!15!black}},
ylabel={closed-loop cost},
axis background/.style={fill=white}
]
\addplot [color=black, dashdotted, line width=2.0pt, mark=*, mark options={solid, black}]
table[row sep=crcr]{%
1353 36.1243007676739\\
1697 36.124300785368\\
2792 36.1243007561747\\
4269 36.1243007649599\\
6177 36.1243007632131\\
7485 36.1243007906212\\
8877 36.1243007777062\\
9335 36.1243007187975\\
};

\addplot [color=black, dashdotted, line width=2.0pt, mark size=5.0pt, mark=x, mark options={solid, black}]
table[row sep=crcr]{%
1353 36.1243007676739\\
1731 26.5006160977074\\
3005 23.5618765913998\\
4085 23.2827444582959\\
6369 23.337151620086\\
8275 23.0612891114376\\
9549 23.0148427334042\\
};

\end{axis}
%
%

\end{tikzpicture}%
}
	\caption{Comparison of the cost of the MPC closed-loop trajectories for different error estimators used for spatial refinement for the time-varying reference $x_{tv}$. $\bullet$ = standard error estimator, {\sf \textbf{x}} = MPC adapted goal oriented error estimator.}
	\label{fig:fvals}
\end{figure}

The advantage of using a goal oriented error estimator that exploits the turnpike structure also shows up when looking at the value of the MPC closed-loop trajetory in dependence of the number of space grid points. It is particularly pronounced for the time varying reference, the respective values can be seen in Figure \ref{fig:fvals} (for other settings see \citep{GrSS20b}). As this figure shows, with the goal oriented error estimator the closed-loop cost is significantly lower for the same number of grid points.

For further numerical examples, including adaptive time discretization as well as unstable and semilinear PDEs we refer to \citep{GrSS19,GrSS20b}.


\section{Topics not Discussed and Open Problems}
\label{sec:open_problems}

\def\N{\mathbb{N}}
\def\R{\mathbb{R}}
\def\U{\mathbb{U}}
\def\X{\mathbb{X}}
\def\Z{\mathbb{Z}}
\def\KK{\mathcal{K}}
\def\NN{\mathcal{N}}
\def\ue{\bar u}
\def\xe{\bar x}
\def\eps{\varepsilon}

In this section we briefly sketch parts of the theory that we did not present in this chapter and discuss some open questions. 
A class of optimal control problems that we did not discuss here are discounted problems, in which the cost is multiplied with a time varying, decreasing factor, typically an exponential of the form $e^{-\delta t}$, $\delta>0$. In discrete time, many of the results discussed in this paper carry over to discounted problems, see \citep{GGHKW18,GMKW20}. In continuous time, however, this still seems to be largely open.

This situation is similar for time-varying problems. Here, one has to distinguish on the one hand between problems with time-invariant data that exhibit a time-varying turnpike, e.g., a periodic orbit instead of an equilibrium. On the other hand, time-varying turnpikes appear naturally in problems depending on time-varying data. Motivated by MPC, both problems have been studied recently, see, e.g., \citep{GrPS18,GruP19,GruM16a,ZaGD17}. However, these investigations are again mostly limited to discrete time.

Both in discrete-time and continuous-time the existing turnpike results focus mostly on the turnpike being in the interior of the constraints, i.e., a comprehensive treatment with active constraints at the turnpike appears to be open. Moreover, problems with mixed-integer inputs and/or hybrid dynamics have only been discussed to a very limited extent, see \citep{tudo:faulwasser20f} for discrete-time dissipativity-based results and \citep{Gugat19a} for  results on infinite-dimensional hyperbolic systems with integrality constraints on the inputs under strict convexity assumptions for the objective.

In recent years the analysis of turnpike properties for infinite-dimensional continuous-time optimal control problems---typically problems governed by PDEs---has received great attention, see, e.g., \citep{FGHS20,Gugat19a,GrSS20,LaTZ20,PorZ13,Zuaz17,Pighin20,porretta2016remarks}. Interestingly, this line of research is usually based on necessary optimality conditions rather than on dissipativity notions with \cite{TreZ18} being one of the few exceptions. A comprehensive theory linking dissipativity and turnpike phenomena for this problem class (as it is available, e.g., for finite-dimensional linear quadratic problems \citep{GruG20}) is still missing, even for the relatively simply linear-quadratic case.

More generally, the connection between the optimality-condition-based and the dissipativity-based approach to turnpike analysis is still not completely understood. Both approaches are usually used in parallel and as an alternative to each other, but many aspects of the relation between these concepts---for instance the question whether there are problem classes for which one of the approaches works while the other fails---have not yet been clarified. Moreover, the implication of turnpike properties on the Hamilton-Jacobi-Bellman-Equation and its solution is not widely explored. First steps have been taken in~\cite{esteve2020turnpike,tudo:faulwasser20a}.

Finally, in mathematical economy there is a huge body of literature on turnpike behavior in stochastic optimal control problems, see \citep{Mari89,KolY12,SSZJ92} for only a small selection of re\-presentative papers. While this property has been established for many specific optimal control problems, a thorough mathematical treatment based, e.g., on stochastic notions of dissipativity, as well as an exploitation of this property for numerical purposes are missing so far.

\vspace*{2mm}

{\bf Acknowledgement:} We thank Manuel Schaller for providing the figures of the numerical examples in Section \ref{sec:Numerics}.

\vspace*{2mm}

\bibliographystyle{agsm}
\bibliography{references_Lars,literature_latin1}
\end{document}